\documentclass[number,citesort,dvips]{arxbj}

\aid{0}
\volume{17}
\issue{1}
\pubyear{2011}
\firstpage{395}
\lastpage{423}
\doi{10.3150/10-BEJ274}

\makeatletter

\newproclaim{Assumption}{Assumption}

\newtheorem{thmm}{Theorem}
\newproclaim{definition}{Definition}
\newtheorem{lemm}{Lemma}
\newtheorem{propm}{Proposition}

\makeatother

\begin{document}
\begin{frontmatter}

\title{Statistical analysis of self-similar conservative fragmentation chains}
\runtitle{Statistical analysis of fragmentation chains}

\begin{aug}
\author[1]{\fnms{Marc} \snm{Hoffmann}\corref{}\thanksref{1}\ead[label=e1]{marc.hoffmann@univ-mlv.fr}} \and
\author[2]{\fnms{Nathalie} \snm{Krell}\thanksref{2}\ead[label=e2]{nathalie.krell@univ-rennes1.fr}}
\runauthor{M. Hoffmann and N. Krell}
\address[1]{Universit\'e Paris Est,
Laboratoire d'Analyse et de Math\'ematiques Appliqu\'ees, CNRS-UMR
8050, 5, Boulevard Descartes, 77454 Marne-la-Vall\'ee Cedex 2,
France.\\
\printead{e1}}
\address[2]{Universit\'e de Rennes 1,
Institut de Recherche math\'ematique de Rennes, CNRS-UMR 6625,
Campus de Beaulieu, B\^atiment 22, 35042 Rennes Cedex, France.\\
\printead{e2}}
\end{aug}

\received{\smonth{3} \syear{2008}}
\revised{\smonth{1} \syear{2010}}

%
\begin{abstract}
We explore statistical inference in self-similar conservative
fragmentation chains when only approximate observations of the sizes of
the fragments below a given threshold are available. This framework,
introduced by Bertoin and Martinez [\textit{Adv. Appl. Probab.}
\textbf{37} (2005) 553--570], is motivated by
mineral crushing in the mining industry. The underlying object that can
be identified from the data is the step distribution of the random walk
associated with a randomly tagged fragment that evolves along the
genealogical tree representation of the fragmentation process. We
compute upper and lower rates of estimation in a parametric framework
and show that in the nonparametric case, the difficulty of the
estimation is comparable to ill-posed linear inverse problems of order
1 in signal denoising.
\end{abstract}

%
\begin{keyword}
\kwd{fragmentation chains}
\kwd{key renewal theorem}
\kwd{nonparametric estimation}
\kwd{parametric}
\end{keyword}

\end{frontmatter}

\section{Introduction}\label{sec1}

\subsection{Motivation}\label{sec11}

Random fragmentation models, commonly used in a variety of physical
models, have their theoretical roots in the works of Kolmogorov
\cite{Kolmogorov} and Filippov \cite{Filippov} (see also
\cite{Aldous,Bertoin1,Krapivsky1,Krapivsky2} and the references therein).
Informally, we imagine an object that falls apart randomly as time
passes. The resulting particles break independently of each other in a
self-similar way. A thorough account on random fragmentation processes
and chains is given in the book by Bertoin \cite{Bertoin1}, a key
reference for this paper.

In this work, we adopt the perspective of statistical inference. We
focus on the quite specific class of \textit{self-similar fragmentation
chains}. The law of a self-similar fragmentation chain is determined by
two components:
\begin{itemize}
\item the dislocation measure, which governs the way that the fragments split;
\item the index of self-similarity, which determines the rate of splitting;
\end{itemize}
see the definition in Section \ref{fragmentation chains}. In this
paper, we postulate a specific observation scheme, motivated by the
mining industry, where the goal is to separate metal from non-valued
components in large mineral blocks by a series of blasting, crushing
and grinding operations.
In this setting, one observes, approximately,
the fragments arising from an initial block of size $m$ only when they
reach a size smaller than some screening threshold, say $\eta>0$; see
\cite{BertoinMartinez} and the references
therein. Asymptotics are taken as the ratio $\varepsilon:= \eta/m$ vanishes.

\subsection{Organization and results of the paper}\label{sec12}

In Section \ref{generaltheory}, we recall the basic tools for the construction
of conservative fragmentation chains, closely following the book by
Bertoin \cite{Bertoin1}. For statistical purposes, our main tool is
the empirical measure ${\mathcal E}_\varepsilon$ of the size of
fragments when they reach a size smaller than a threshold
$\varepsilon$ in the limit $\varepsilon\rightarrow0$. We highlight
the fact that ${\mathcal E}_\varepsilon$ captures information about
the dislocation measure through the L\'evy measure $\pi$ of a
randomly tagged fragment associated with the fragmentation process.

In Section \ref{sec3}, we give a rate of convergence for the empirical
measure ${\mathcal E}_\varepsilon$ toward its limit in Theorem
\ref{rate empiricalmeasure}, extending former results (under more
stringent assumptions) of Bertoin and Martinez
\cite{BertoinMartinez}. The rate is of the form $\varepsilon
^{1/2-\ell(\pi)}$, where $\ell(\pi)>0$ can be made arbitrarily
small under suitable 
exponential moment conditions for $\pi$. We additionally consider %
the more realistic framework of observations with limited accuracy,
where each fragment is actually known up to a systematic stochastic
error of order $\sigma\ll\varepsilon$. We construct estimators
related to
functionals of $\pi$ in the absolutely continuous case. In the
parametric case (Theorem \ref{lowerparametric}), we establish that the
best achievable rate is $\varepsilon^{1/2}$, in the particular case of
binary fragmentations, where a particle splits into two blocks at each
step exactly. We construct a convergent estimator in a general setting
(Theorem \ref{upperparametric}) with an error of order $\varepsilon
^{1/2-\ell'(\pi)}$ for another $\ell'(\pi)>0$ that can be made
arbitrarily small under appropriate assumptions on the density of $\pi
$ near $0$ and $+\infty$.
In the nonparametric case, we construct an estimator that achieves
(Theorem \ref{upperbound}) a rate of the form $(\varepsilon^{1-\ell
''(\pi)})^{s/(2s+3)}$, where $s>0$ is the local smoothness of the
density of $\pi$, up to appropriate rescaling. Except for the factor
$\ell''(\pi) >0$, we obtain the same rate as for ill-posed inverse
problems of degree 1.


\section{Statistical model} \label{generaltheory}

\subsection{Fragmentation chains} \label{fragmentation chains}

A fragmentation chain can be constructed as follows. We start with a
state space
\[
{\mathcal S}^\downarrow:= \Biggl\{\mathbf{s} = (s_1,s_2,\ldots), s_1 \geq s_2
\geq\cdots\geq0, \sum_{i=1}^\infty s_i \leq1 \Biggr\}.
\]
A point $\mathbf{s} \in{\mathcal S}^\downarrow$ is interpreted as the
collection of (decreasing) sizes of fragments originating from a single
(unit) mass. We also specify the following two quantities:
\begin{itemize}
\item a finite dislocation measure $\nu$, that is, a finite measure
$\nu(\mathrm{d}\mathbf{s})$ on ${\mathcal S}^\downarrow$;
\item a parameter of self-similarity, $\alpha\geq0$.
\end{itemize}
A fragmentation chain with parameter of self-similarity $\alpha$ and
dislocation measure $\nu$ is a Markov process $X= (X(t),t \geq0 )$
with value in ${\mathcal S}^\downarrow$. Its evolution can be
described as follows:
a fragment with size $x$ lives for an exponential time with parameter
$x^\alpha\nu({\mathcal S}^\downarrow)$ and then splits and gives
rise to a family of smaller fragments distributed as $x\xi$, where
$\xi$ is distributed according to $\nu(\cdot)/\nu({\mathcal
S}^\downarrow)$. We denote by $\mathbb{P}_m$ the law of $X$ started
from the
initial configuration $(m,0,\ldots)$ with $m \in(0,1]$. Under
$\mathbb{P}
_m$, the law of $X$ is entirely determined by $\alpha$ and $\nu(\cdot
)$; see Theorem 3 of Bertoin \cite{Bertoin2}. To ensure that
everything is well defined, 
the following assumptions on the dislocation measure $\nu(\mathrm{d}\mathbf{s})$
of $X$ are in force throughout the paper.
\renewcommand{\theAssumption}{\Alph{Assumption}}
\begin{Assumption}\label{assuA}
We have $\nu({\mathcal S}^\downarrow)=1$
and $\nu(s_1\in(0,1) )=1$.
\end{Assumption}

In our setting, Assumption \ref{assuA} is standard; see Bertoin \cite{Bertoin1}.
We will repeatedly use the representation of fragmentation chains as
random infinite marked trees. Let
\[
{\mathcal U} := \bigcup_{n=0}^\infty\mathbb{N}^n
\]
denote the infinite genealogical tree (with $\mathbb{N}^0:=\{\varnothing
\}$)
associated with $X$ as follows: to each node $u \in{\mathcal U}$, we
set a mark
\[
(\xi_u, a_u, \zeta_u),
\]
where $\xi_u$ is the size of the fragment labeled by $u$, $a_u$ is its
birth-time and $\zeta_u$ is its life-time.
We have the following identity between point measures on $(0,+\infty)$:
\[
\sum_{i=1}^\infty1_{\{X_i(t)>0 \}}\delta_{X_i(t)} = \sum_{u \in
{\mathcal U}}1_{ \{t \in[a_u,a_u+\zeta_u) \}}\delta_{\xi_u},\qquad t\geq0,
\]
with $X(t)= (X_1(t),X_2(t),\ldots)$ and where $\delta_x$ denotes the
Dirac mass at $x$.
Finally, $X$ has the following branching property: for every fragment
$\mathbf{s} = (s_1,\ldots) \in{\mathcal S}^\downarrow$ and every $t
\geq0$, the distribution of $X(t)$ given $X(0)=\mathbf{s}$ is the same as
the decreasing rearrangement of the terms of independent random
sequences $X^{(1)}(t)$, $X^{(2)}(t),\ldots,$ where, for each $i$,
$X^{(i)}(t)$ is distributed as $X(t)$ under~$\mathbb{P}_{s_i}$.

\subsection{Observation scheme} \label{observation scheme}
Keeping in mind the motivation of mineral crushing, we consider the
fragmentation under
$\mathbb{P}:=\mathbb{P}_1$, initiated with a unique block of size
$m=1$, and we
observe the process stopped at the time when all the fragments become
smaller than some given threshold $\varepsilon>0$, so we have data
$\xi_u$, for every $u \in{\mathcal U}_{ \varepsilon}$,
with
\[
{\mathcal U}_\varepsilon:= \{u \in{\mathcal U}, \xi_{u-}\geq
\varepsilon, \xi_u <\varepsilon\},
\]
where we denote by $u-$ the parent of the fragment labeled by $u$.
We will further assume that the total mass of the fragments remains
constant through time, as follows.
\begin{Assumption}[(Conservative property)]\label{assuB}
We have
$\nu(\sum_{i=1}^\infty s_i=1 )=1$.
\end{Assumption}

We next consider a test function $g(\cdot)$ integrated against the
empirical measure
\[
{\mathcal E}_\varepsilon(g):=\sum_{u \in{\mathcal U_\varepsilon
}}\xi_u g(\xi_u/\varepsilon).
\]
Indeed, under Assumption \ref{assuB}, we have
%
\begin{equation} \label{conservative}
\sum_{u \in{\mathcal U}_\varepsilon}\xi_u=1 \qquad\mathbb{P}\mbox
{-almost surely,}
\end{equation}
so ${\mathcal E}_\varepsilon(g)$ appears as a weighted empirical
version of $g(\cdot)$. Note that the empirical measure ${\mathcal
E}_\varepsilon$ depends only on the size of the fragmentation and is
thus independent of the self-similarity parameter $\alpha$. Bertoin
and Martinez show in
\cite{BertoinMartinez}, Corollary 1, that under mild assumptions on
$\nu(\cdot)$, the random variable
${\mathcal E}_\varepsilon(g)$ converges to
\[
{\mathcal E}(g):=\frac{1}{c(\nu)}\int_{0}^1 \frac{g(a)}{a} \int
_{{\mathcal S}^\downarrow} \sum_{i=1}^\infty s_i 1_{\{s_i <a\}}\nu
(\mathrm{d}\mathbf{s})\,\mathrm{d}a
\]
in $L^1(\mathbb{P})$ as $\varepsilon\rightarrow0$, with
$c(\nu)=-\int_{{\mathcal S}^\downarrow}\sum_{i=1}^\infty s_i\log
s_i \nu(\mathrm{d}\mathbf{s})$, tacitly assumed to be well defined.
This suggests a strategy for recovering information about $\nu(\cdot
)$ by choosing suitable test functions $g(\cdot)$. In Section \ref
{rate of convergence L2}, we will show that the convergence also holds
in $L^2(\mathbb{P})$ and we will exhibit a rate of convergence, which
is a
crucial issue if statistical results are sought.

\subsection{First estimates} \label{first estimates}
From now on, we assume that we have data
%
\begin{equation} \label{obsscheme}
X_\varepsilon:= (\xi_u, u \in{\mathcal U}_\varepsilon)
\end{equation}
and we specialize in the estimation of $\nu(\cdot)$. Clearly, the
data give no information about the parameter of self-similarity 
$\alpha$ that we consider as a nuisance parameter.
Assumptions \ref{assuA} and \ref{assuB} are in force. At this stage, we can relate
${\mathcal E}(g)$ to a more appropriate quantity by means of the
so-called \textit{tagged fragment} approach.

\textit{The randomly tagged fragment}. Let us first consider the
homogenous case $\alpha=0$. Assume that we can ``tag'' a point at
random according to a uniform distribution on the initial fragment and
imagine that we can follow the evolution of the fragment that contains
this point. Let us denote by $ (\chi(t), t \geq0 )$ the process of
the size of the fragment that contains the randomly chosen point.
This fragment is a typical observation in our data set $X_\varepsilon$
and it appears at time
\[
T_\varepsilon:=\inf\{t \geq0, \chi(t) < \varepsilon\}.
\]
Bertoin \cite{Bertoin1} shows that the process
$\zeta(t):=-{\log}\chi(t)$
is a subordinator with L\'evy measure
%
\begin{equation} \label{defLevy}
\pi(\mathrm{d}x):=\mathrm{e}^{-x}\sum_{i=1}^\infty\nu(-{\log s_i} \in \mathrm{d}x).
\end{equation}
We can anticipate that the information we get from $X_\varepsilon$ is
actually information about the L\'evy measure $\pi(\mathrm{d}x)$ of $\zeta(t)$
obtained via 
$\zeta(T_\varepsilon)$. The dislocation measure
$\nu(\mathrm{d}\mathbf{s})$ and $\pi(\mathrm{d}x)$ are related by (\ref{defLevy}), which reads
%
\begin{equation} \label{characterizationpi}
\int_{{\mathcal S}^\downarrow} \sum_{i=1}^\infty s_i f(s_i)\nu
(\mathrm{d}\mathbf{s})=\int_{(0,+\infty)}f(\mathrm{e}^{-x})\pi(\mathrm{d}x)
\end{equation}
for any suitable $f(\cdot)\dvtx [0,1]\rightarrow[0,+\infty)$. In
particular, by Assumption \ref{assuB} and the fact that $\nu({\mathcal
S}^\downarrow)=1$, $\pi(\mathrm{d}x)$ is a probability measure, hence $\zeta
(t)$ is a compound Poisson process. Informally, a~typical observation
takes the form $\zeta(T_\varepsilon)$, which is the value of a
subordinator with L\'evy measure $\pi(\mathrm{d}x)$ at its first passage time
strictly above $-{\log\varepsilon}$.
The case $\alpha\neq0$ is a bit more involved and reduces to the
homogenous case by a time change; see Bertoin
\cite{Bertoin2,Bertoin1}. In terms of the limit of the empirical measure ${\mathcal
E}_\varepsilon(g)$, we equivalently have
\[
{\mathcal E}(g) = \frac{1}{c(\pi)}\int_0^1 \frac{g(a)}{a} \pi
(-{\log a},+\infty) \,\mathrm{d}a = \frac{1}{c(\pi)}\int_0^{+\infty} g(\mathrm{e}^{-x})
\pi(x,+\infty) \,\mathrm{d}x
\]
with $c(\pi)=\int_{(0,+\infty)}x \pi(\mathrm{d}x)$. The representation of
${\mathcal E}(g)$ as an integral with respect to $\pi$ will prove
technically convenient. Except in the binary case (a particular case of
interest, see Section~\ref{DiscussionBinarycase}), knowledge of $\pi
(\cdot)$ does not, in general, allow us to recover $\nu(\cdot)$.

\textit{Measurements with limited accuracy}. It is unrealistic to assume
that we can observe exactly the sizes $\xi_u$ of the fragments. This
becomes even more striking if the dislocation splits at a given time
into infinitely many fragments of non-zero size, a situation that we do
not discard in principle.
Therefore, we replace (\ref{obsscheme}) by the more realistic
observation scheme
$X_{\varepsilon,\sigma}:= (\xi_u^{(\sigma)}, u \in{\mathcal
U}_{\varepsilon, \sigma} )$
with
\[
{\mathcal U}_{\varepsilon, \sigma}:= \bigl\{u \in{\mathcal U}, \xi
_{u-}^{(\sigma)} \geq\varepsilon, \xi_u^{(\sigma)} < \varepsilon\bigr\}
\]
and
%
\begin{equation} \label{defobs}
\xi_u^{(\sigma)}:=\xi_u +\sigma U_u.
\end{equation}
The random variables $(U_u, u \in{\mathcal U})$ are identically
distributed and account for a systematic experimental microstructure
noise in the measurement of $X_\varepsilon$,
independent of $X_\varepsilon$. We assume, furthermore, that for every
$u \in{\mathcal U}$,
\[
|U_u| \leq1 \quad\mbox{and}\quad \mathbb{E}[U_u]=0.
\]
The noise level $0 \leq\sigma= \sigma(\varepsilon) \ll\varepsilon
$ is assumed to be known and represents the accuracy level of the
statistician. The observations $\xi_u+\sigma U_u$ are further
discarded below a threshold $\sigma\leq t_\varepsilon\leq\varepsilon
$, beyond which they become irrelevant, leading to the modified
empirical measure
\[
{\mathcal E}_{\varepsilon,\sigma}(g):=\sum_{u \in{\mathcal U}_{
\varepsilon, \sigma}} 1_{\{\xi_u^{(\sigma)} \geq t_\varepsilon\}}
\xi_u^{(\sigma)} g \bigl(\xi_u^{(\sigma)}/\varepsilon\bigr).
\]
In the sequel, we take $t_\varepsilon= \gamma_0\varepsilon$ for some
(arbitrary) $0 < \gamma_0 < 1$ and assume further that
$\sigma\leq\frac{1}{2}t_\varepsilon$.

\section{Main results}\label{sec3}

\subsection{A rate of convergence for the empirical measure}
\label{rate of convergence L2}

\begin{definition}
For $\kappa>0$, we say that a non-lattice probability measure
$\pi(\mathrm{d}x)$ defined on $[0,+\infty)$ belong to ${\Pi}(\kappa)$ if
$
\int_{[0,+\infty)}\mathrm{e}^{\kappa x} \pi(\mathrm{d}x) < +\infty$. We
set $\Pi(\infty):=\bigcap_{\kappa>0}\Pi(\kappa)$.
\end{definition}

For $m >0$, let
\[
{\mathcal C}(m):= \Bigl\{g\dvtx[0,1]\rightarrow\mathbb{R}\mbox{, continuous, }
\|g\|
_\infty:=\sup_x|g(x)| \leq m \Bigr\}
\]
and
%
\[
{\mathcal C}'(m):= \Bigl\{g \in{\mathcal C}(m)\dvtx[0,1]\rightarrow\mathbb
{R}\mbox{, differentiable, }\|g'\|_\infty:=\sup_x|g'(x)| \leq m \Bigr\}.
\]
Our first result exhibits explicit rates in the convergence ${\mathcal
E}_\varepsilon(g) \rightarrow{\mathcal E}(g)$ as $\varepsilon
\rightarrow0$, extending Bertoin \cite{Bertoin1}, Proposition 1.12.
\begin{thmm} \label{rate empiricalmeasure}
We work under 
Assumptions \textup{\ref{assuA}} and \textup{\ref{assuB}}.
Let $1 < \kappa\leq\infty$ and assume that $\pi
\in\Pi(\kappa)$.
\begin{itemize}
\item For every $m>0$ and $1 \leq\mu<
\kappa$, we have
%
\begin{equation} \label{convergence}
\sup_{g \in{\mathcal C}(m)} \mathbb{E}\bigl[ \bigl({\mathcal E}_\varepsilon(g)
-{\mathcal E}(g) \bigr)^2 \bigr]=\mathrm{o} \bigl(\varepsilon^{\mu/(\mu+1)} \bigr).
\end{equation}
\item The convergence (\ref{convergence}) remains valid if we replace
${\mathcal E}_{\varepsilon}(\cdot)$ by ${\mathcal E}_{\varepsilon,
\sigma}(\cdot)$ and ${\mathcal C}(m)$ by ${\mathcal C}'(m)$. The
following additional error term must then be incorporated: for any $0 <
\mu< \kappa$, we have
%
\begin{equation} \label{convergencebis}
\sup_{g \in{\mathcal C}'(m)} \mathbb{E}\bigl[ \bigl({\mathcal E}_{\varepsilon
,\sigma
}(g) -{\mathcal E}_\varepsilon(g) \bigr)^2 \bigr]
=\mathrm{o} (\varepsilon^{\mu/2} )+{\mathcal O}(\sigma\varepsilon^{-1}).
\end{equation}
\end{itemize}
\end{thmm}

\subsection{Statistical estimation} \label{statpreliminaries}

We study the estimation of $\pi(\cdot)$
by constructing estimators based on ${\mathcal E}_\varepsilon(\cdot)$
or, rather, ${\mathcal E}_{\varepsilon, \sigma}(\cdot)$.
We need the following regularity assumption.
\begin{Assumption}\label{assuC}
The probability $\pi(\mathrm{d}x)$ is absolutely
continuous with respect to the Lebesgue measure: $\pi(\mathrm{d}x)=\pi(x)\,\mathrm{d}x$.
Moreover, its density function $x \leadsto\pi(x)$ is continuous on
$(0,+\infty)$ and satisfies
$\limsup_{x \rightarrow+\infty}\mathrm{e}^{\vartheta x}\pi(x)<+\infty$ for
some $\vartheta\geq1$.
\end{Assumption}

We distinguish two cases: the \textit{parametric case}, where we estimate
a linear functional of $\pi(\cdot)$ of the form
\[
m_k(\pi):=\int_0^{+\infty}x^k\pi(x)\,\mathrm{d}x,\qquad k=1,2,\ldots,
\]
and the
\textit{nonparametric case}, where we estimate the function $x \leadsto
\pi(x)$ pointwise. In the latter case, it will prove convenient to
assess the local smoothness properties of $\pi(\cdot)$ on a
logarithmic scale. Henceforth, we consider the mapping
%
\begin{equation} \label{logtransformation}
a \leadsto\beta(a):=a^{-1}\pi(-{\log a}),\qquad a \in(0,1).
\end{equation}
In the nonparametric case, we estimate $\beta(a)$ for every $a \in(0,1)$.

\subsection{The parametric case}\label{sec33}

\textit{Preliminaries}. For $k \geq1$, we estimate
\[
m_k(\pi):=\int_0^{+\infty}x^k \pi(x)\,\mathrm{d}x = \int_0^1 \log(1/a)^k
\beta(a)\,\mathrm{d}a
\]
by the correspondence (\ref{logtransformation}), implicitly assumed to
be well defined. We first focus on the case $k=1$. Choose a
sufficiently smooth test function $f(\cdot)\dvtx[0,1]\rightarrow\mathbb
{R}$ such
that $f(1)=0$ and let $g(a):=-af'(a)$. Clearly,
%
\begin{eqnarray}\label{energieg}
{\mathcal E}(g) &=& \frac{1}{c(\pi)}\int_0^1 \frac{g(a)}{a} \pi
(-{\log a},+\infty) \,\mathrm{d}a \nonumber\\[-8pt]\\[-8pt]
&=& -\frac{1}{m_1(\pi)}\int_0^1 f'(a)\int_0^a \beta(u) \,\mathrm{d}u\, \mathrm{d}a
=\frac{1}{m_1(\pi)}\int_0^1 f(a)\beta(a)\,\mathrm{d}a.\nonumber
\end{eqnarray}
%
Formally, taking $f(\cdot)\equiv1$ would identify $1/m_1(\pi)$ since
$\beta(\cdot)$ integrates to one, but this choice is forbidden by the
boundary condition $f(1)=0$. We shall instead consider a family of
regular functions that are close to the constant function $1$ while
satisfying $f(1)=0$.

\textit{Construction of the approximating functions}.
Let $f_{\gamma}\dvtx[0,1]\rightarrow\mathbb{R}$ with $0<\gamma<1$ be a family
of smooth functions satisfying the following conditions:
\begin{itemize}
\item we have $f_{\gamma}(a)=1$ for $a \leq1-\gamma$ and $f_{\gamma}(1)=0$;
\item we have
%
\begin{equation} \label{bornitudederiveesf}
\sup_{\gamma>0} (\|f_\gamma\|_\infty+\gamma\|f'_{\gamma}\|_\infty
+\gamma^2 \|f_{\gamma}''\|_\infty) <+\infty;
\end{equation}
\item for every $k \geq1$ and some $\delta>0$, we have
%
\begin{eqnarray} \label{controlefgammaorigine}
&&\sup_{\gamma>0}\sup_{a \in(0,1)} \biggl\{\gamma^2|{\log a}|^k \bigl(a^{-1}
|f_\gamma(1-a)|+|f'_\gamma(1-a)| \bigr)+\biggl(\frac{\gamma}{a}\biggr)^{1+\delta
}f_\gamma(1-a) \biggr\}\qquad\nonumber\\[-8pt]\\[-8pt]
&&\quad<+\infty.\nonumber
\end{eqnarray}
\end{itemize}
The family $(f_\gamma,\gamma>0)$ mimics the behaviour of the target
function $f_0(a)=1$ for $0 \leq a < 1$ and $f_0(1)=0$ as $\gamma
\rightarrow0$. Condition (\ref{controlefgammaorigine}) is technical
(and probably not optimal). An explicit choice of a family $(f_\gamma
,\gamma>0)$ satisfying (\ref{bornitudederiveesf}) and (\ref
{controlefgammaorigine}) is given by
\[
f_\gamma(a):=
\cases{
1 &\quad if $a \leq1-\gamma$, \cr
10 \biggl(\dfrac{1-a}{\gamma} \biggr)^3-15 \biggl(\dfrac{1-a}{\gamma} \biggr)^4+6 \biggl(\dfrac
{1-a}{\gamma} \biggr)^5 &\quad if $1-\gamma\leq a < 1$,
\cr
0 &\quad if $a=1$,}
\]
but other choices are obviously possible.

\textit{Construction of an estimator}.
We are now ready to give an estimator of the first moment $m_1(\pi)$
of $\pi$ and, more generally, of any moment $m_k(\pi)$, $k \geq1$.
For a parametrization $\gamma:=\gamma_\varepsilon\rightarrow0$ to
be specified later, we set
\[
g_{\gamma_\varepsilon}(a):=-a f_{\gamma_\varepsilon}'(a),\qquad a \in(0,1).
\]
By Theorem \ref{rate empiricalmeasure}, we expect ${\mathcal
E}_{\sigma,\varepsilon}(g_{\gamma_{\varepsilon}})$ to be close to
${\mathcal E}(g_{\gamma_{\varepsilon}})$ which, in turn, is equal to
$m_1(\pi)^{-1}\int_0^1 f_{\gamma_\varepsilon}(a)\beta(a)\,\mathrm{d}a$, by
(\ref{energieg}). Since $f_{\gamma_\varepsilon} \approx1$ and
$\beta(\cdot)$ is a density function, by appropriate regularity
assumptions on $\pi$, we may further expect this last quantity to be
close to $1/m_1(\pi)$. We therefore set
%
\begin{equation} \label{estimationmoyenne}
\widehat m_{1,\varepsilon}:=\frac{1}{{\mathcal E}_{\varepsilon
,\sigma} (g_{\gamma_\varepsilon} )}
\end{equation}
for an estimator of $m_1(\pi)$.
More generally, for $k>1$, we define successive moment estimators as
follows. Set $h_{\gamma_\varepsilon}(a):=f_{\gamma_\varepsilon
}(1-a)\log(1/a)^k$ and $\widetilde g_{\gamma_\varepsilon}(a):=-a
h_{\gamma_\varepsilon}'(a)$. The same heuristics as before lead to
the estimator
\[
\widehat m_{k,\varepsilon}:=\frac{{\mathcal E}_{\varepsilon,\sigma}
(\widetilde g_{\gamma_\varepsilon} )}{{\mathcal E}_{\varepsilon
,\sigma} (g_{\gamma_\varepsilon} )}.
\]

\textit{Upper rates of convergence}.
We can describe the performance
of $\widehat m_{k,\varepsilon}$ under an additional decay condition on
$\pi(\cdot)$ near the origin.
\begin{definition}
For $\kappa> 0$, we say
that the probability $\pi(\cdot)$ belong to the class ${\mathcal
R}(\kappa)$ if
$\limsup_{x \rightarrow0}x^{-\kappa+1}\pi(x) < +\infty$.
We set
${\mathcal R}(\infty):=\bigcap_{\kappa>0}{\mathcal R}(\kappa)$.
\end{definition}

We obtain the following upper bound, under more stringent regularity
assumptions on $\pi$ than in Theorem \ref{rate empiricalmeasure}.
\begin{thmm} \label{upperparametric}
We work under 
Assumptions \textup{\ref{assuA}}, \textup{\ref{assuB}} and \textup{\ref{assuC}}.
\begin{itemize}
\item For the estimation of $m_1(\pi)$, assume $\kappa_1 \geq4$ and
$\kappa_2 > 1$.
\item For the estimation of $m_k(\pi)$ with $k \geq2$, assume $\kappa
_1 \geq4$ and $\kappa_1 > \kappa_2 > 1$.
\end{itemize}
For any $1 \leq\mu< \kappa_1$, let $\widehat m_{k,\varepsilon}$ be
specified by $\gamma_\varepsilon:=\varepsilon^{\mu/(\mu+1)(2\kappa
_2+1)}$. The family
\[
\bigl(\varepsilon^{-\mu/(\mu+1)} \bigr)^{\kappa_2/(2\kappa_2+1)} \bigl(\widehat
m_{k,\varepsilon}-m_k(\pi) \bigr)
\]
is tight provided that
\[
\pi\in\Pi(\kappa_1) \cap{\mathcal R}(\kappa_2)
\]
and $\sigma\varepsilon^{-3}$ remains bounded.
\end{thmm}

Some remarks: The convergence of $\widehat m_{k,\varepsilon}$ to
$m_{k}(\pi)$ is of course no surprise, by (\ref{convergence}).
However, the dependence on $\varepsilon$ in the test function
$g_{\gamma_\varepsilon}(\cdot)$ (in particular, $g_{\gamma
_\varepsilon}(\cdot)$ 
is unbounded as $\varepsilon\rightarrow0$) requires a slight
improvement of Theorem \ref{rate empiricalmeasure}. This can be done
thanks to Assumption \ref{assuC}; see Proposition \ref
{reinforcemenempiricalmeasure} in Section \ref{proof upperparametric}.
The requirement $\sigma\varepsilon^{-3} = {\mathcal O}(1)$ ensures
that the additional term coming from the approximation of ${\mathcal
E}_\varepsilon(\cdot)$ by ${\mathcal E}_{\sigma,\varepsilon}(\cdot
)$ is negligible.

\textit{Lower rates of convergence}.
Our next result shows that the exponent
\[
\frac{\mu}{\mu+1} \frac{\kappa_2}{2\kappa_2+1} \leq\frac{1}{2}
\]
in the rate of convergence of Theorem \ref{upperparametric} is nearly
optimal, to within an arbitrarily small polynomial order.
\begin{definition} \label{lowerbound}
Let $\pi_0(\cdot)$ satisfy the assumptions of Theorem \ref{upperparametric}.
The rate $0 < v_\varepsilon\rightarrow0$ is a lower rate of convergence
for estimating $m_k(\pi_0)$ if there exists a family $\pi_\varepsilon
(\cdot)$ satisfying the assumptions of Theorem \ref{upperparametric}
and a constant $c>0$ such that
%
\begin{equation} \label{deflowerbound}
\liminf_{\varepsilon\rightarrow0}\inf_{F_\varepsilon}\max_{\pi
\in\{\pi_0,\pi_\varepsilon\}} \mathbb{P}[v_\varepsilon
^{-1}|F_\varepsilon-m_k(\pi)| \geq c ] >0,
\end{equation}
where the infimum is taken (for every $\varepsilon$) over all
estimators constructed from $X_{\varepsilon,\sigma}$ at level
$\varepsilon$. 
\end{definition}

Definition \ref{lowerbound} expresses a kind of local min--max
information bound: given $\pi_0(\cdot)$, one can find $\pi
_\varepsilon(\cdot)$ such that no estimator can discriminate between
$\pi_0(\cdot)$ and $\pi_\varepsilon(\cdot)$ at a rate faster than
$v_\varepsilon$. We further restrict our attention to binary
fragmentations; see Section \ref{DiscussionBinarycase}.
In that case, the dislocation measure satisfies $\nu(s_1+s_2\neq1)=0$
and, because of the conservation Assumption \ref{assuB}, can be represented as
%
\begin{equation} \label{rho-representation}
\nu(\mathrm{d}\mathbf{s})=\rho(\mathrm{d}s_1)\delta_{1-s_1}(\mathrm{d}s_2),
\end{equation}
where $\rho(\cdot)$ is a probability measure on $[1/2,1]$.
\begin{Assumption}[(Binary case)]\label{assuD}
The probability measure $\rho
(\cdot)$ associated with $\pi(\cdot)$ is absolutely continuous and
its density function is bounded away from zero.
\end{Assumption}
\begin{thmm} \label{lowerparametric} Assume that the fragmentation is
binary and work under Assumption \textup{\ref{assuD}}. In the same setting as in Theorem
\ref{upperparametric}, the rate $\varepsilon^{1/2}$ is a lower rate
of convergence for estimating $m_k(\pi)$.
\end{thmm}


\subsection{The nonparametric case}\label{sec34}

\textit{Preliminaries}.
Under local smoothness assumptions on the function 
$\beta(\cdot)$, we estimate
$\beta(a)$ for every $a\in(0,1)$. Given $s>0$, we say that $\beta
(\cdot)$ belongs to the H\"older class $\Sigma(s)$ if there exists a
constant $c>0$ such that
\[
\bigl|\beta^{(n)}(y)-\beta^{(n)}(x) \bigr|\leq c|y-x|^{\{s\}}
\]
with $s=n+\{s\}$, where $n$ is a non-negative integer and $\{s\} \in(0,1]$.
We also need to relate $\beta(\cdot)$ to the decay of its
corresponding L\'evy measure $\pi(\cdot)$.
Again abusing notation, we identify $\Pi(\kappa)$ with the set of
$\beta(\cdot)$ such that
$\mathrm{e}^x\beta(\mathrm{e}^{-x})\,\mathrm{d}x \in\Pi(\kappa)$,
thanks to the inverse of (\ref{logtransformation}), and likewise for
${\mathcal R}(\kappa)$.

\textit{Construction of an estimator}. We construct an estimator of $\beta
(\cdot)$ in the same way as for the parametric case: for $a \in(0,1)$
and a normalizing factor $0 < \gamma_\varepsilon\rightarrow0$, set
\[
\varphi_{\gamma_\varepsilon,a}(x):=\gamma_\varepsilon^{-1}\varphi
\bigl((x-a)/\gamma_\varepsilon\bigr),
\]
where $\varphi$ is a smooth function with support in $(0,1)$ that
satisfies the following oscillating property: for some integer $N \geq1$,
%
\begin{equation} \label{cancellationproperty}
\int_0^1\varphi(a)\,\mathrm{d}a = 1,\qquad \int_0^1 a^k \varphi(a)\,\mathrm{d}a=0,\qquad k=1,\ldots, N.
\end{equation}
The function $\varphi_{\gamma_\varepsilon,a}$ thus plays the role of
a kernel centred around $a$. Set
\[
h_{a,\varepsilon}(x) = -x \varphi'_{\gamma_\varepsilon, a}(x),\qquad x
\in(0,1).
\]
We have
\[
{\mathcal E}(h_{a,\varepsilon}) = \frac{1}{m_1(\pi)}\int_0^1\varphi
_{\gamma_\varepsilon,a}(x)\beta(x)\,\mathrm{d}x,
\]
by (\ref{energieg}). By letting $h_\varepsilon\rightarrow0$ with an
appropriate rate as $\varepsilon\rightarrow0$, we expect this term to
be close to $\beta(a)/m_1(\pi)$. Eventually, we can eliminate the
denominator by means of our preliminary estimator $\widehat
m_{1,\varepsilon}$. Our nonparametric estimator of $\beta(a)$ thus
takes the form
\[
\widehat\beta_\varepsilon(a):=\widehat m_{1,\varepsilon} {\mathcal
E}_{\varepsilon,\sigma} (h_{a,\varepsilon} ),\qquad a\in(0,1),
\]
where $\widehat m_{1,\varepsilon}$ is the estimator of $m_1(\pi)$
defined in (\ref{estimationmoyenne}).

\textit{Upper rates of convergence}. We have the following
result.
\begin{thmm} \label{upperbound}
We work under
Assumptions \textup{\ref{assuA}}, \textup{\ref{assuB}}
and \textup{\ref{assuC}}. Let $\kappa_1 \geq4$ and $\kappa_2 >1$. For
any $1 \leq\mu< \kappa_1$,
let $\widehat\beta_\varepsilon(\cdot)$ be
specified by $\gamma_\varepsilon:=\varepsilon^{\mu/(\mu
+1)(2s+3)}$. For every $a \in(0,1)$, the family
\[
\bigl(\varepsilon^{-\mu/(\mu+1)} \bigr)^{s/(2s+3)} \bigl(\widehat\beta
_\varepsilon(a)-\beta(a) \bigr)
\]
is tight, provided that
\[
\beta\in\Sigma(s) \cap\Pi(\kappa_1) \cap{\mathcal R}(\kappa_2)
\]
for $0 < s < \min\{N,3\kappa_2\}$ and $\sigma\varepsilon^{-3}$
remains bounded.
\end{thmm}

A proof of the (near) optimality, in the sense of the lower bound
Definition \ref{lowerbound} and in the spirit of Theorem \ref
{lowerparametric}, is presumably a delicate problem that lies beyond
the scope of the paper; see 
Section~\ref{discussion nonparametric}.

\section{Discussion} \label{Discussion}

\subsection{Binary fragmentations} \label{DiscussionBinarycase}

The case of binary fragmentations is the simplest, yet is an important
model of random fragmentation, where a particle splits into two blocks
at each step (see, e.g., \cite{Filippov,BrennanDurett}).
By using representation (\ref{rho-representation}), if we further
assume that $\rho(\mathrm{d}s_1)=\rho(s_1)\,\mathrm{d}s_1$ is absolutely continuous, then
so is $\pi(\mathrm{d}x)=\pi(x)\,\mathrm{d}x$ and we have
%
\begin{equation} \label{rhopi-correspondence}
\pi(x)=\mathrm{e}^{-2x} \bigl(\rho(\mathrm{e}^{-x})1_{[0,\log2]}(x)+\rho(1-\mathrm{e}^{-x})1_{(\log
2, +\infty)}(x) \bigr)
\end{equation}
for $x\in[0,+\infty)$
and
%
\begin{equation} \label{rhobetta-correspondence}\beta(a)=a \bigl(\rho
(a)1_{[1/2,1]}(a)+\rho(1-a)1_{[0, 1/2)}(a) \bigr),\qquad a \in[0,1].
\end{equation}
In particular, the regularity properties of $\beta(\cdot)$ are
obtained from the local smoothness of $\rho(\cdot)$ and its behaviour
near $1$. For instance, if $\rho(1-a) ={\mathcal O}(a^{\kappa-1})$
near the origin, for some $\kappa>0$, then
\[
\pi\in\Pi(\kappa)\cap{\mathcal R}(\kappa).
\]

\subsection[Concerning Theorem 1]{Concerning Theorem \protect\ref{rate empiricalmeasure}}\label{sec42}

Theorem \ref{rate empiricalmeasure} readily extends to error
measurements of the form $\mathbb{E}[ |{\mathcal E}_\varepsilon(g)
-{\mathcal
E}(g) |^p ]$ with $1 \leq p \leq2$. The rate becomes $\varepsilon
^{-\mu p/2(\mu+1)}$ in (\ref{convergence}) and $\sigma^p\varepsilon
^{-p}$ in (\ref{convergencebis}) under the less stringent condition
$\mu< \kappa/2p$.

Generally speaking, in (\ref{convergence}), we obtain the (normalized)
rate $\varepsilon^{\mu/2(\mu+1)}$ for any $\mu<\kappa$.
Intuitively, we have a number of observations that should be of order
$\varepsilon^{-1}$, so the expected rate would rather be $\varepsilon
^{1/2}$. Why can we not obtain the rate $\varepsilon^{1/2}$, or simply
$\varepsilon^{\kappa/2(\kappa+1)}$? The proof in Section~\ref{proof
rate empiricalmeasure} shows that we lose quite a lot of information %
when applying Sgibnev's result (see Proposition \ref
{renouvellementSgibnev} in Section \ref{Sgibnevth}) on the key renewal theorem for
a random walk with step distribution $\pi(\cdot)$ in the limit $\log
(1/\varepsilon)\rightarrow+\infty$.

Proposition \ref{renouvellementSgibnev} ensures that if $\pi(\cdot)$
has exponential moments up to order $\kappa$, then we can guarantee in
the renewal theorem the rate $\mathrm{o}(\varepsilon^\mu)$ for any $\mu<
\kappa$ with some uniformity in the test function, a crucial point for
the subsequent statistical applications. It is presumably possible to
improve this rate to ${\mathcal O}(\varepsilon^{\kappa})$ by using %
Ney's result \cite{Ney}. However, a careful glance at the proof of
Theorem \ref{rate empiricalmeasure} shows that we would then lose an
extra logarithmic term when replacing $\varepsilon^{\mu/2(\mu+1)}$
by $\varepsilon^{\kappa/(2\kappa+1)}$. More generally, exhibiting
exact rates of convergence in Theorem \ref{rate empiricalmeasure}
remains a delicate issue: the key renewal theorem is sensitive to a
modification of the distribution outside a neighbourhood of $+\infty$;
see, for example, Asmussen \cite{Asmussen}, page 196.

\subsection[Concerning Theorems 2 and 4]{Concerning Theorems \protect\ref{upperparametric} and
\protect\ref{upperbound}} \label{discussion nonparametric}

In the parametric case, we obtain the rate
\[
\bigl(\varepsilon^{\mu/(\mu+1)} \bigr)^{\kappa_2/(2\kappa_2+1)} \qquad\mbox{for
all } \mu< \kappa_1,
\]
%
which can be made arbitrary close to the lower bound $\varepsilon
^{1/2}$ by assuming $\kappa_1$ and $\kappa_2$ to be large enough. The
factor $\mu/(\mu+1)$ comes from Theorem \ref{rate empiricalmeasure},
whereas the factor $\kappa_2/(2\kappa_2+1)$ arises when using the
technical assumption $\pi\in{\mathcal R}(\kappa_2)$. We do not know
how to improve this.

In the nonparametric case, the situation is a bit different than in the
parametric case: we now obtain the rate
%
\begin{equation} \label{resultoptimal}
\bigl(\varepsilon^{\mu/(\mu+1)} \bigr)^{s/(2s+3)} \qquad\mbox{for all } \mu< \kappa_1
\end{equation}
for the estimation of $\beta(a)$ for any $a \in(0,1)$. In the limit
$\kappa_1 \rightarrow+\infty$, it becomes $\varepsilon^{s/(2s+3)}$,
which can be related to more classical models in the nonparametric
literature. Informally, a function of $d$ variables with degree of
smoothness $s$ observed in noise under the action of a smoothing
operator of degree $\nu$ (e.g., $\nu$-fold integration) can be
recovered with optimal rate
$\varepsilon^{s/(2s+2\nu+d)}$; see, for instance, \cite{Tsybakov2}.
Here, we have $d=1$ and $\nu=1$ by the representation (\ref
{energieg}), so formula~(\ref{resultoptimal}) is consistent with the
general nonparametric theory. This advocates in favour of the (near)
optimality of the result in the sense of Definition \ref{lowerbound},
but a complete proof lies beyond the scope of the paper.
\subsection{The Crump--Mode--Jagers alternative}\label{sec44}

As suggested by a referee, the statistical problem can be reformulated
alternatively in terms of the Crump--Mode--Jagers (CMJ) branching
process. Consider a transformed point process $(\tau_1, \tau_2,\ldots
)$ defined by $\tau_i=-{\log s_i}$ for $\mathbf{s} = (s_1,s_2,\ldots) \in
{\mathcal S}^\downarrow$. The sequence $(\tau_1, \tau_2,\ldots)$
describes the consecutive ages at childbearing for the individual
assumed to be born at time zero. In our setting, the resulting CMJ
process is supercritical with Malthusian parameter $1$ since $\mathrm{e}^{-\tau
_1}+\mathrm{e}^{-\tau_2}+\cdots= 1$.

Let $\sigma_u = -{\log\xi_u}$. We may now interpret $\sigma_u$ as the
individual forming the coming generation at time $t =
-{\log\varepsilon}
$. The empirical measure ${\mathcal E}_\varepsilon$ now has the representation
\begin{eqnarray*}
{\mathcal E}_\varepsilon(g)& = & \sum_{u \in{\mathcal U},\sigma
_u-\tau_u \leq t < \sigma_u}\mathrm{e}^{-\sigma_u}g(\mathrm{e}^{-\sigma_u+t})\\
& = & \mathrm{e}^{-t}\sum_{u \in{\mathcal U},\sigma_u-\tau_u \leq t < \sigma
_u}\mathrm{e}^{-\sigma_u+t}g(\mathrm{e}^{-\sigma_u+t})
\end{eqnarray*}
and the last sum can be expressed in terms of a population size with
random characteristics; see \cite{JN}. This yields
another interpretation of our statistical approach in terms of
branching processes, presumably more useful in other settings.

\section{Proofs}\label{sec5}

We will repeatedly use the convenient notation $a_\varepsilon\lesssim
b_\varepsilon$ if $0 < a_\varepsilon\leq c b_\varepsilon$ for some
constant $c>0$ which may depend on $\pi(\cdot)$ and on the constant
$m$ appearing in the definition of the class ${\mathcal C}(m)$ or
${\mathcal C}'(m)$. Any other dependence on other ancillary quantities
will be obvious from the context.
A function $g \in{\mathcal C}(m)$ is tacitly defined on the whole real
line by setting $g(a)=0$ for $a \notin[0,1]$.

\subsection{Preliminaries: Rates of convergence in the key renewal
theorem} \label{Sgibnevth}

We state a special case of Sgibnev's result \cite{Sgibnev} 
on uniform rates of convergence in the key renewal theorem, an
essential tool for this paper. Let $F(\mathrm{d}x)$ be a non-lattice probability
distribution with positive mean $m$ and renewal function $\mathbb
{F}=\sum_{n=0}^\infty F^{n\star}$ with $F^{0\star}:=\delta_0$,
$F^{1\star}:=F$ and $F^{(n+1)\star}:=F\star F^{n\star}$, $n\geq0$.
We denote by
$T(F)$ the $\sigma$-finite measure with density function
\[
\int_{(x,+\infty)}F(\mathrm{d}u) 1_{[0,+\infty)}(x) -\int_{(-\infty
,x]}F(\mathrm{d}u) 1_{(-\infty,0)}(x)
\]
and define $T^2(F):=T (T(F) )$. Let $\varphi(\cdot)\dvtx \mathbb
{R}\rightarrow
[0,+\infty)$ be a submultiplicative function, that is, such that
$\varphi(0)=1, \varphi(x+y)\leq\varphi(x)\varphi(y)$.
We then have (see, e.g., \cite{Phillips}, Section 6)
\begin{eqnarray*}
-\infty<r_1:\!&=&\lim_{x\rightarrow-\infty}\frac{\log\varphi
(x)}{x}\\
&\leq&\lim_{x\rightarrow+\infty}\frac{\log\varphi
(x)}{x}=:r_2<+\infty.
\end{eqnarray*}
\begin{Assumption}\label{assuE}
We have $r_1 \leq0 \leq r_2$ and there exists
$r\dvtx\mathbb{R}\rightarrow\mathbb{R}$, an integrable function such
that the following
conditions are fulfilled:
\[
\sup_x|r(x)|\varphi(x)<+\infty,\qquad \lim_{|x|\rightarrow\infty
}r(x)\varphi(x)=0,
\]
\[
\lim_{x\rightarrow+\infty}\varphi(x)\int_{[x, +\infty)}r(u)\,\mathrm{d}u
=\lim_{x\rightarrow-\infty}\varphi(x)\int_{(-\infty, x]}r(u)\,\mathrm{d}u=0
\]
and
$\int_{\mathbb{R}} \varphi(x)T^2(F)(\mathrm{d}x)<\infty$.
We call $\varphi(\cdot)$ a \textit{rate function} and $r(\cdot)$ a
\textit{dominating function}.
\end{Assumption}

Sgibnev's result takes the following form.
\begin{propm}[(\cite{Sgibnev}, Theorem 5.1)]\label
{renouvellementSgibnev} We work under 
Assumption \textup{\ref{assuE}}. Then
\[
\lim_{|t|\rightarrow\infty}\varphi(t)\sup_{\psi,|\psi(x)|\leq
|r(x)|} \biggl|\psi\star\mathbb{F}(t)-m^{-1}\int_{\mathbb{R}}\psi(x)\,\mathrm{d}x \biggr|=0.
\]
\end{propm}

\subsection[Proof of Theorem 1]{Proof of Theorem \protect\ref{rate empiricalmeasure}}
\label{proof rate empiricalmeasure}

\textit{Step} 1: \textit{A preliminary decomposition}.
We first use the fact that for $\eta> \varepsilon$, during the
fragmentation process, the unobserved state $X_\eta$ necessarily
anticipates the state $X_\varepsilon$. The choice $\eta= \eta
(\varepsilon)$ will follow later. This yields the following representation:
\[
{\mathcal E}_\varepsilon(g) = \sum_{v \in{\mathcal U}_\eta}\xi
_v\sum_{w \in{\mathcal U}}1_{ \{\xi_{v}\tilde\xi_{w-}^{(v)}\geq
\varepsilon, \xi_{v}\tilde\xi_{w}^{(v)} < \varepsilon\}}\tilde\xi
_{w}^{(v)}g \bigl(\xi_v\tilde\xi_w^{(v)}/\varepsilon\bigr),
\]
where, for each label $v \in{\mathcal U}_\eta$ and conditional on
$X_\eta$, a new independent fragmentation chain $(\tilde\xi
^{(v)}_{w}, w\in{\mathcal U})$ is started, thanks to the branching
property; see Section \ref{fragmentation chains}. Now, define
\[
\lambda_\eta(v):=1_{\{\xi_{v-}\geq\eta,\xi_v <\eta\}}\xi_v
\]
and
\[
Y_\varepsilon(v,g):=\sum_{w \in{\mathcal U}} 1_{ \{\xi_{v}\tilde
\xi_{w-}^{(v)}\geq\varepsilon, \xi_{v}\tilde\xi_{w}^{(v)} <
\varepsilon\}}\tilde\xi_{w}^{(v)}g \bigl(\xi_v\tilde\xi
_w^{(v)}/\varepsilon\bigr).
\]
We obtain the decomposition of ${\mathcal E}_\varepsilon(g)-{\mathcal
E}(g)$ as the sum of a centred and a bias term:
\[
{\mathcal E}_\varepsilon(g)-{\mathcal E}(g) = M_{\varepsilon,\eta
}(g)+B_{\varepsilon, \eta}
\]
with
\[
M_{\varepsilon,\eta}(g):=\sum_{v \in{\mathcal U}} \lambda_\eta
(v) \bigl(Y_\varepsilon(v,g)-\mathbb{E}[Y_\varepsilon(v,g) | \lambda
_\eta(v) ] \bigr)
\]
and
\[
B_{\varepsilon,\eta}(g):=\sum_{v \in{\mathcal U}} \lambda_\eta
(v) \bigl(\mathbb{E}[Y_\varepsilon(v,g) | \lambda_\eta(v) ]-{\mathcal
E}(g) \bigr),
\]
where we have used the conservative property (\ref{conservative}) in
order to incorporate the limit term ${\mathcal E}(g)$ into 
the sum in $v$.

\textit{Step} 2: \textit{The term $M_{\varepsilon,\eta}(g)$}.
Conditional on the $\sigma$-field generated by the random variables $
(1_{\{\xi_{v-}\geq\eta\}}\xi_v,v-\in{\mathcal U} )$, the variables
$(Y_\varepsilon(v,g), v \in{\mathcal U})$ are independent. Therefore,
%
\begin{equation} \label{inegmg}
\mathbb{E}[M_{\varepsilon,\eta}(g)^2 ]\leq\sum_{v \in{\mathcal
U}}\mathbb{E}
[\lambda_\eta(v)^2\mathbb{E}[Y_\varepsilon(v,g)^2 | \lambda_\eta
(v)] ].
\end{equation}
Thus, we first need to control the conditional variance of
$Y_\varepsilon(v,g)^2$ given $\lambda_\eta(v) = u$, for $0 \leq u
\leq\eta$, since $\mathbb{P}$-almost surely, $\lambda_\eta(v) \leq
\eta
$. Moreover, we have $Y_{\varepsilon}(v,g)=0$ on the event $\{\lambda
_\eta(v) < \varepsilon\}$, hence we may assume that $\varepsilon\leq
u \leq\eta$.

To this end, we will use the following representation property.
\begin{lemm} \label{representation}
Let $f(\cdot)\dvtx [0,+\infty)\rightarrow[0,+\infty)$. Then
%
\begin{equation}
\mathbb{E}\biggl[\sum_{v \in{\mathcal U}_\eta}\xi_v f(\xi_v) \biggr]= \mathbb
{E}^\star[f
(\chi(T_\eta) ) ],
\end{equation}
where $\chi(t)=\exp(-\zeta(t) )$ and $ (\zeta(t), t \geq0 )$ is a
subordinator with L\'evy measure $\pi(\cdot)$ defined on an
appropriate probability space $(\Omega^\star, \mathbb{P}^\star)$ and
\[
T_\eta:=\inf\{t\geq0, \zeta(t) > -{\log\eta}\}.
\]
\end{lemm}

The proof readily follows the construction of the randomly tagged
fragment as elaborated in the book by Bertoin \cite{BertoinLevy} and
is thus omitted.
We plan to bound the right-hand side of (\ref{inegmg}) using 
Lemma \ref{representation}. For $0 < \varepsilon\leq u \leq\eta$,
we have
\begin{eqnarray*}
\mathbb{E}[Y_\varepsilon(v,g)^2 | \lambda_\eta(v)=u ]
&=& \mathbb{E}\biggl[ \biggl(\sum_{w \in{\mathcal U}_{\varepsilon/u}}\widetilde
\xi
_w^{(v)}g \bigl(\varepsilon u^{-1}\widetilde\xi_w^{(v)} \bigr) \biggr)^2 \Big| \lambda
_\eta(v)=u \biggr] \\
&\leq&\mathbb{E}\biggl[\sum_{w \in{\mathcal U}_{\varepsilon/u}}\widetilde
\xi
_w^{(v)}g \bigl(\varepsilon u^{-1}\widetilde\xi_w^{(v)} \bigr)^2 \Big| \lambda
_\eta(v)=u \biggr],
\end{eqnarray*}
where we have used Jensen's inequality combined with (\ref
{conservative}). Applying Lemma \ref{representation}, we derive
%
\begin{equation} \label{passagepi}
\mathbb{E}[Y_\varepsilon(v,g)^2 | \lambda_\eta(v)=u ] \leq\mathbb
{E}^\star\bigl[g
\bigl(u\varepsilon^{-1}\mathrm{e}^{-\zeta(T_{\varepsilon/
u})} \bigr)^2 \bigr].
\end{equation}
Let $U(\cdot)$ denote the renewal function associated with the
subordinator $(\zeta(t), t \geq0)$. By \cite{BertoinLevy},
Proposition 2, Chapter III, the right-hand side of (\ref{passagepi})
is equal to
\begin{eqnarray*}
&&\int_{ [0, -{\log(\varepsilon/u)} )}\mathrm{d}U(s)\int_{(-{\log(\varepsilon
/u)} - s,+\infty)}g (u \varepsilon^{-1}\mathrm{e}^{-x-s} )^2\pi(\mathrm{d}x) \\
&&\quad= \int_{ [0,-{\log(\varepsilon/u) })}\mathrm{d}U(s) \int_{{\mathcal
S}^\downarrow}\sum_{i=1}^\infty s_i 1_{\{s_i < \varepsilon u^{-1}
\mathrm{e}^s\}} g (s_iu \varepsilon^{-1}\mathrm{e}^{-s} )^2\nu(\mathrm{d}\mathbf{s}) \\
&&\quad\lesssim\frac{1}{c(\pi)}\|g\|_\infty^2\log(u/\varepsilon),
\end{eqnarray*}
where we have successively used the representation
(\ref{characterizationpi}) and the upper bound $U(s) \lesssim
s/c(\pi)$; see, for instance, \cite{BertoinLevy}, Proposition 1,
Chapter III.
Therefore, for $\varepsilon\leq u \leq\eta$,
\[
\mathbb{E}[Y_\varepsilon(v,g)^2 | \lambda_\eta(v)=u ]\lesssim\frac
{1}{c(\pi)}\|g\|_\infty^2\log(\eta/\varepsilon).
\]
Going back to (\ref{inegmg}), since $\lambda_\eta(v)^2\leq\eta
\lambda_\eta(v)$ and again using (\ref{conservative}),
we readily derive
%
\begin{equation} \label{ordremg}
\mathbb{E}[M_{\varepsilon,\eta}(g)^2 ] \lesssim\frac{1}{c(\pi)}\|
g\|
_\infty^2 \eta\log(\eta/\varepsilon) \lesssim\eta\log(\eta
/\varepsilon).
\end{equation}

\textit{Step} 3: \textit{The bias term $B_{\varepsilon, \eta}(g)$}.
First, note that
\[
\mathbb{E}[Y_\varepsilon(v,g) | \lambda_\eta(v) ]=\xi_v^{-1}
\mathbb{E}_{ \xi_v}
[{\mathcal E}_\varepsilon(g) ],
\]
$\mathbb{P}$-almost surely, so
%
\begin{equation} \label{decompositionbiais}
B_{\varepsilon,\eta}(g)=\sum_{v \in{\mathcal U}}\lambda_{\eta
}(v) \bigl(\xi_v^{-1} \mathbb{E}_{\xi_v} [{\mathcal E}_\varepsilon(g)
]-{\mathcal E}(g) \bigr).
\end{equation}
Conditioning on the mark of the parent $v-=\omega$ of $v$ and applying
the branching property, we get that $\mathbb{E}_{ \xi_v} [{\mathcal
E}_\varepsilon(g) ]$ can be written as
\[
\mathbb{E}_{ \xi_v} \Biggl[\sum_{\omega\in{\mathcal U}}1_{ \{\widehat
{\xi
}_{\omega} \geq\varepsilon\}}\widehat{\xi}_\omega\int_{{\mathcal
S}^\downarrow}\sum_{i=1}^\infty1_{ \{\widehat{\xi}_{\omega} s_i
<\varepsilon\}}s_i g (\widehat{\xi}_\omega s_i \varepsilon^{-1}
)\nu(\mathrm{d}\mathbf{s}) \Biggr],
\]
where the $(\widehat{\xi}_w, w \in{\mathcal U})$ are the sizes of
the marked fragments
of a fragmentation chain with same dislocation measure $\nu(\cdot)$,
independent of $(\xi_v, v \in{\mathcal U})$. Set
\[
H_g(z):=\int_{{\mathcal S}^\downarrow}\sum_{i=1}^\infty1_{ \{s_i <
\mathrm{e}^{-z} \}}s_i g (s_i \mathrm{e}^z )\nu(\mathrm{d}\mathbf{s}),\qquad z\geq0.
\]
It follows that $\mathbb{E}_{ \xi_v} [{\mathcal E}_\varepsilon(g) ]$ is
equal to
\begin{eqnarray*}
&&\mathbb{E}_{\xi_v} \Biggl[\sum_{n=0}^\infty\sum_{|\omega|=n}1_{ \{\log
\widehat{\xi}_\omega\geq\log\varepsilon\}}\widehat{\xi}_\omega
H_g(\log\widehat{\xi}_\omega-\log\varepsilon) \Biggr] \\
&&\quad= \xi_v \mathbb{E}\Biggl[\sum_{n=0}^\infty\sum_{|\omega|=n}1_{ \{\log
\widehat{\xi}_\omega\geq\log(\varepsilon/\rho) \}} \widehat{\xi
}_\omega H_g \bigl(\log\widehat{\xi}_\omega-\log(\varepsilon/\rho) \bigr)
\Biggr]_{\rho= \xi_v},
\end{eqnarray*}
by self-similarity, with the notation $|\omega|=n$ if $\omega=(\omega
_1,\ldots,\omega_n) \in{\mathcal U}$. Using \cite{Bertoin1},
Proposition 1.6, we finally obtain
\[
\mathbb{E}_{ \xi_v} [{\mathcal E}_\varepsilon(g) ]=\xi_v\sum
_{n=0}^\infty
\mathbb{E}\bigl[1_{ \{S_n \leq\log(\rho/\varepsilon) \}}H_g \bigl(\log(\rho
/\varepsilon)-S_n \bigr) \bigr]_{\rho=\xi_v},
\]
where $S_n$ is a random walk with step distribution $\pi(\mathrm{d}x)$. Note
that this can also be written as
%
\begin{equation} \label{eqrepresentation}
\xi_v^{-1} \mathbb{E}_{ \xi_v} [{\mathcal E}_\varepsilon(g) ] =
\mathbb
{F}\star\psi\bigl(\log(\xi_v/\varepsilon) \bigr),
\end{equation}
where $\mathbb{F} = \sum_{n=0}^\infty\pi^{n\star}$ denotes the
renewal measure associated with the probability measure $\pi$ and
$\psi(z) = 1_{z \leq0}H_g(-z)$. In order to bound
\[
\xi_v^{-1} \mathbb{E}_{\xi_v} [{\mathcal E}_\varepsilon(g)
]-{\mathcal E}(g),
\]
we plan to apply a version of the renewal theorem with explicit rate of
convergence as given in Sgibnev \cite{Sgibnev}; see Proposition \ref
{renouvellementSgibnev} in Section \ref{Sgibnevth}. We take a rate
function $\varphi(z):=\exp(\mu' z)$ for some arbitrary $\mu' <
\kappa/2$, a dominating function $r(z):= \mathrm{e}^{-\kappa|z|}$ and set
$F:=\pi$ in Proposition \ref{renouvellementSgibnev}. We can write,
for $z < 0$,
\[
H_g(-z) = 1_{ \{z\leq0 \}}\int_{(-z,+\infty)}g (\mathrm{e}^{-x-z} )\pi(\mathrm{d}x),
\]
by (\ref{characterizationpi}).
Since $g(\cdot)$ has support in $[0,1]$ and $\pi\in\Pi(\kappa)$,
\[
|H_g(-z) |\leq\int_{(-z,+\infty)} |g (\mathrm{e}^{-x-z} ) |\pi(\mathrm{d}x) \lesssim
\mathrm{e}^{\kappa z}.
\]
Therefore, $|1_{\{z\leq0\}}H_g(-z)|\lesssim r(z)$ for all $z \in
\mathbb{R}$.
Since $\kappa> 2\mu'$, Assumption \ref{assuE} of Proposition \ref
{renouvellementSgibnev} is readily checked. Now, let $A>0$ (depending
only on $\kappa$, $m$ and $\pi(\cdot))$ such that, if
$\log(\xi_v/\varepsilon) \geq A$, then, by Proposition \ref
{renouvellementSgibnev},
%
\begin{equation} \label{renouvellement}
\biggl|\xi_v^{-1} \mathbb{E}_{ \xi_v} [{\mathcal E}_\varepsilon(g)
]-\frac{1}{\mathbb{E}
^\star[S_1 ]}\int_0^{+\infty}H_g(z)\,\mathrm{d}z \biggr| \leq\biggl(\frac{\varepsilon
}{\xi_v} \biggr)^{\mu'}.
\end{equation}
We next note that
\[
\frac{1}{\mathbb{E}^\star[S_1 ]}\int_0^{+\infty}H_g(z)\,\mathrm{d}z =
{\mathcal E}(g).
\]
Introducing the family of events $ \{\log(\xi_v/\varepsilon) \geq A
\}$ in the sum (\ref{decompositionbiais}), we obtain the following
decomposition:
\[
B_{\varepsilon,\eta}(g)^2 \lesssim I + \mathit{II}
\]
with
\[
I:=\sum_{v \in{\mathcal U}_\eta}\xi_v1_{ \{\log(\xi
_v/\varepsilon) > A \}} \bigl(\xi_v^{-1} \mathbb{E}_{ \xi_v} [{\mathcal
E}_\varepsilon(g) ]-{\mathcal E}(g) \bigr)^2
\]
and
\[
\mathit{II}:=\sum_{v \in{\mathcal U}_\eta}\xi_v1_{ \{\log(\xi
_v/\varepsilon) \leq A \}} \bigl(\xi_v^{-1} \mathbb{E}_{ \xi_v}
[{\mathcal
E}_\varepsilon(g) ]-{\mathcal E}(g) \bigr)^2.
\]
By (\ref{renouvellement}), we have
\[
I \leq\varepsilon^{2\mu'}\sum_{v \in{\mathcal U}_\eta}1_{ \{
-{\log\xi_v} < -A+\log(1/\varepsilon) \}}\xi_{v}\exp\bigl(2\mu'(-{\log
\xi_v}) \bigr).
\]
Integrating with respect to $\mathbb{P}$ and applying Lemma \ref
{representation}, in the same way as in step 2, we have
\begin{eqnarray*}
\mathbb{E}[I ] & \leq & \varepsilon^{2\mu'} \mathbb{E}^\star\bigl[\mathrm{e}^{2\mu
' \zeta
(T_\eta)} \bigr]\\
& = & \varepsilon^{2\mu'}\int_{[0,-{\log\eta})}\mathrm{d}U(s)\int_{(-{\log\eta
-s},+\infty)}\mathrm{e}^{2\mu'(s+x)}\pi(\mathrm{d}x) \\
& \leq & \varepsilon^{2\mu'} \int_{[0,-{\log\eta})} \mathrm{e}^{2\mu' s}\,\mathrm{d}U(s)
\lesssim(\varepsilon\eta^{-1} )^{2\mu'} \log(1/\eta)
\end{eqnarray*}
for small enough $\varepsilon$ and where we have used $\pi\in\Pi
(\kappa)$ with $2\mu' < \kappa$. For the term $\mathit{II}$, we first note
that by (\ref{conservative})
and self-similarity,
\[
\mathbb{E}_{\xi_v} \biggl[\sum_{u \in{\mathcal U}_\varepsilon}\widehat
\xi_u
\biggr]=\xi_v,\qquad \mathbb{P}_{\xi_v}\mbox{-almost surely,}
\]
hence
%
\begin{equation} \label{point important}
\bigl(\xi_v^{-1}\mathbb{E}_{ \xi_v} [{\mathcal E}_\varepsilon(g)
]-{\mathcal
E}(g) \bigr)^2\leq4\|g\|_\infty^2,\qquad \mathbb{P}_{\xi_v}\mbox{-almost surely.}
\end{equation}
In the same way as for the term $I$, we derive
\begin{eqnarray*}
\mathbb{E}[\mathit{II} ] & \lesssim & \mathbb{E}\biggl[\sum_{v \in{\mathcal U}_\eta
}\xi_v 1_{ \{
-{\log\xi_v} \geq-A+\log(1/\varepsilon) \}} \biggr]\\
&=& \mathbb{P}^\star[\zeta(T_\eta) \geq-A+\log(1/\varepsilon) ]\\
&\leq&\int_{[0,-{\log\eta})}\mathrm{d}U(s) \int_{(-A+\log(1/\varepsilon
)-s,+\infty)}\pi(\mathrm{d}x)\\
&\lesssim&\varepsilon^{\mu'} \log(1/\eta)
\end{eqnarray*}
for small enough $\varepsilon$. Using 
all of the estimates together, we conclude that
%
\begin{equation} \label{conclusionbiais}
\mathbb{E}[B_{\varepsilon,\eta}(g)^2 ] \lesssim
\bigl(\varepsilon^{\mu'}+(\varepsilon\eta^{-1})^{2\mu'}
\bigr)\log(1/\eta).
\end{equation}

\textit{Step} 4: \textit{Proof of} (\ref{convergence}).
Using the estimates (\ref{ordremg}) and (\ref{conclusionbiais}), we have
\begin{eqnarray*}
\mathbb{E}\bigl[ \bigl({\mathcal E}_\varepsilon(g)-{\mathcal E}(g) \bigr)^2 \bigr]
&\lesssim& \mathbb{E}
[M_{\varepsilon,\eta}(g)^2 ]+\mathbb{E}[B_{\varepsilon,\eta}(g)^2
]\\
& \lesssim & \eta\log(\eta/\varepsilon)+ (\varepsilon\eta^{-1}
)^{2\mu'} \log(1/\eta)+\varepsilon^{\mu'} \log(1/\eta).
\end{eqnarray*}
%
The choice $\eta(\varepsilon):=\varepsilon^{2\mu'/(2\mu'+1)}$
yields the rate
\[
\varepsilon^{\min\{2\mu'/(2\mu'+1),\mu'\}}\log(1/\varepsilon)
\qquad\mbox{for any } 0 < \mu' < \kappa/2.
\]
We thus obtain a rate of the form $\mathrm{o}(\varepsilon^{\mu/(\mu+1)})$ for
any $1 \leq\mu< \kappa$. The conclusion follows.

\textit{Step} 5: \textit{Proof of} (\ref{convergencebis}). We plan to use the
following decomposition:
\[
{\mathcal E}_{\varepsilon,\sigma}(g)-{\mathcal E}_{\varepsilon}(g) =
I + \mathit{II}
\]
with
\[
I:= \sum_{u \in{\mathcal U}} \bigl(1_{\{\xi_{u-}^{(\sigma)} \geq
\varepsilon, \xi_u^{(\sigma)} < \varepsilon\}}-1_{\{ \xi_{u-}\geq
\varepsilon, \xi_u< \varepsilon\}} \bigr)\widetilde\xi_{u}^{(\sigma)}
g \bigl(\xi_u^{(\sigma)}/\varepsilon\bigr)
\]
and
\[
\mathit{II}:=\sum_{u \in{\mathcal U}_\varepsilon} \bigl(\widetilde\xi_u^{(\sigma
)} g \bigl(\xi_u^{(\sigma)}/\varepsilon\bigr)- \xi_u g (\xi_u/\varepsilon) \bigr),
\]
where we have set $\widetilde\xi_u^{(\sigma)}:=\xi_u^{(\sigma
)}1_{\{\xi_u^{(\sigma)} \geq t_\varepsilon\}}$. Clearly,
\begin{eqnarray*}
\bigl|1_{\{\xi_{u-}^{(\sigma)} \geq\varepsilon, \xi_u^{(\sigma)} <
\varepsilon\}}-1_{\{ \xi_{u-}\geq\varepsilon, \xi_u< \varepsilon\}
} \bigr| &\leq& 1_{\{\xi_{u-}^{(\sigma)} \geq\varepsilon, \xi_{u-}<
\varepsilon\}}+ 1_{\{\xi_{u}^{(\sigma)} < \varepsilon, \xi_{u}
\geq\varepsilon\}}\\
&&{} +  1_{\{\xi_{u-} \geq\varepsilon, \xi_{u-}^{(\sigma)}<
\varepsilon\}} +1_{\{\xi_{u} < \varepsilon, \xi_{u}^{(\sigma)}
\geq\varepsilon\}}.
\end{eqnarray*}
Let $\delta> \sigma/\varepsilon$ and $\omega=u$ or $u-$. Since
$|U_\omega|\leq1$ for every $\omega$, we can readily check that
\[
\bigl\{ \xi_\omega^{(\sigma)}\geq\varepsilon, \xi_\omega<
\varepsilon\bigr\} \subset\{(1-\delta)\varepsilon\leq\xi_\omega
<\varepsilon\}
\]
and
\[
\bigl\{ \xi_\omega\geq\varepsilon, \xi_\omega^{(\sigma)} <
\varepsilon\bigr\} \subset\{ \varepsilon\leq\xi_\omega<(1+\delta
)\varepsilon\}.
\]
It follows that $|I| \leq \mathit{III} + \mathit{IV}$
with
\[
\mathit{III}:=\sum_{u \in{\mathcal U}}1_{\{(1-\delta)\varepsilon\leq\xi
_{u-} \leq\varepsilon(1+\delta)\}} \bigl|\widetilde\xi_{u}^{(\sigma)}
g \bigl(\xi_u^{(\sigma)}/\varepsilon\bigr) \bigr|
\]
and
\[
\mathit{IV}:=\sum_{u \in{\mathcal U}}1_{\{ (1-\delta)\varepsilon\leq\xi
_{u} \leq(1+\delta)\varepsilon\}} \bigl|\widetilde\xi_{u}^{(\sigma)} g
\bigl(\xi_u^{(\sigma)}/\varepsilon\bigr) \bigr|.
\]
By choosing $\delta$ to be small enough, we may (and will) assume that
$\widetilde\xi_u^{(\sigma)} \lesssim\xi_u$. Conditioning on the
mark of the parent $u-=v$ of $u$, using the branching property,
Jensen's inequality and the conservative Assumption \ref
{conservative}, we conclude that
$\mathbb{E}[\mathit{III}^2 ] $ is less than
\begin{eqnarray*}
&& \mathbb{E}\Biggl[\sum_{v \in{\mathcal U}}1_{\{(1-\delta)\varepsilon
\leq\xi
_{v} \leq\varepsilon(1+\delta)\}}\xi_{v} \int_{{\mathcal
S}^\downarrow} \sum_{i=1}^\infty s_i g \bigl(\varepsilon^{-1}(\xi_v
s_i+\sigma U_v) \bigr)^2\nu(\mathrm{d}\mathbf{s}) \Biggr] \\
&&\quad= \mathbb{E}\biggl[\sum_{\omega\in{\mathcal U}}1_{\{(1-\delta
)\varepsilon\leq
\xi_{\omega} \leq\varepsilon(1+\delta)\}}\xi_{\omega} G_1(\xi
_\omega) \biggr]
\end{eqnarray*}
with
\[
G_1(a):=\int_{{\mathcal S}^\downarrow}\sum_{i=1}^\infty s_i \mathbb{E}\bigl[g
\bigl(\varepsilon^{-1}(a s_i+\sigma U) \bigr)^2 \bigr]\nu(\mathrm{d}\mathbf{s})
\]
and $U$ distributed as the $U_\omega$.
Likewise,
\[
\mathbb{E}[\mathit{IV}^2 ] \lesssim\mathbb{E}\biggl[\sum_{u \in{\mathcal U}}1_{\{
(1-\delta
)\varepsilon\leq\xi_{u} \leq\varepsilon(1+\delta)\}}\xi_{u}
G_2(\xi_u) \biggr]
\]
with
$G_2(a):=\mathbb{E}[g (\varepsilon^{-1}(a+\sigma U) )^2 ]$.
For $i=1,2$, the crude bound $|G_i(a)| \leq\|g\|_\infty^2$ and the
genealogical representation argument used in step 3 enable us 
to bound either $\mathbb{E}[\mathit{III}^2]$ or $\mathbb{E}[\mathit{IV}^2]$ by
\[
\|g\|_\infty^2\sum_{n=0}^\infty{\mathbb P}^\star[-{\log(1+\delta)}
\leq S_n - \log(1/\varepsilon) \leq-{\log(1-\delta)} ],
\]
where $S_n$ is a random walk with step distribution $\pi(\cdot)$. We
proceed as in step 3 and apply Proposition \ref
{renouvellementSgibnev}. The above term converges to
\[
m_1(\pi)^{-1}\log\biggl(\frac{1+\delta}{1-\delta} \biggr) \lesssim\delta
\]
uniformly in $\delta$, provided that $\delta$ is bounded, at rate
$\varepsilon^{\mu'}$ for any $0 < \mu' < \kappa/2$, and is thus of
order $\delta+ \varepsilon^{\mu'}$. We next turn to the term $\mathit{II}$.
We have $\mathit{II}:=V+\mathit{VI}+\mathit{VII}$ with
\begin{eqnarray*}
V&:=&\sum_{u \in{\mathcal U}_{ \varepsilon}} \xi_u \bigl(g \bigl(\xi
_u^{(\sigma)}/\varepsilon\bigr)-g (\xi_u/\varepsilon) \bigr),\\
\mathit{VI}&:=&\sigma\sum_{u \in{\mathcal U}_{ \varepsilon}} U_u 1_{\{\xi
_u^{(\sigma)} \geq t_\varepsilon\}} g \bigl(\xi_u^{(\sigma)}/\varepsilon
\bigr),\\
\mathit{VII}&:=&-\sum_{u \in{\mathcal U}_{ \varepsilon}} \xi_u 1_{\{\xi
_u^{(\sigma)} < t_\varepsilon\}} g \bigl(\xi_u^{(\sigma)}/\varepsilon\bigr).
\end{eqnarray*}
From $g\in{\mathcal C}'(m)$, (\ref{conservative}), Jensen's
inequality and a Taylor expansion, 
we derive that
\[
\mathbb{E}[V^2] \leq\|g'\|_\infty^2\sigma^2\varepsilon^{-2}.
\]
From $|U_u|\leq1$ and the inclusion $\{\xi_u^{(\sigma)} \geq
t_\varepsilon\} \subset\{\xi_u \geq t_\varepsilon-\sigma\}$, we derive
\[
\mathbb{E}[\mathit{VI}^2 ] \leq\|g\|_\infty^2 \frac{\sigma
^2}{(t_\varepsilon
-\sigma)^2}\mathbb{E}\biggl[ \biggl(\sum_{u \in{\mathcal U}_\varepsilon}\xi_u
\biggr)^2 \biggr]
\lesssim\frac{\sigma^2}{\varepsilon^2},
\]
where we have used the fact that $t_\varepsilon=\gamma_0\varepsilon$
with $0<\gamma_0<1$ and $\sigma\leq t_\varepsilon/2$. Likewise, the
inclusion $\{\xi_u^{(\sigma)} < t_\varepsilon\} \subset\{\xi_u
\leq t_\varepsilon+\sigma\}$ and Lemma \ref{representation} yield
\[
\mathbb{E}[\mathit{VII}^2 ] \leq\|g\|_\infty^2\mathbb{P}^\star[
-{\log\chi}
(T_\varepsilon
)>-{\log(t_\varepsilon+\sigma) }] \lesssim\varepsilon^{\mu'} \log
(1/\varepsilon)
\]
for any $0 < \mu' < \kappa/2$, along the same lines as for the bound
of the right-hand side of (\ref{passagepi}) in step~2. Putting all of
the estimates together with, for instance, $\delta:=\sigma
/2\varepsilon$, we finally obtain a rate of the form
\[
\varepsilon^{\mu'}\log(1/\varepsilon)+\sigma\varepsilon^{-1}
\qquad\mbox{for any } 0 < \mu' < \kappa/2,
\]
which can be written as $\mathrm{o}(\varepsilon^{\mu/2})+{\mathcal O}(\sigma
\varepsilon^{-1})$ for any $0 < \mu< \kappa$. We thus obtain (\ref
{convergencebis}) and the proof of Theorem \ref{rate empiricalmeasure}
is complete.

\subsection[Proof of Theorem 2]{Proof of Theorem \protect\ref{upperparametric}}
\label{proof upperparametric}

\textit{Preliminaries}. We begin with a technical lemma.
\begin{lemm} \label{betaborne}
We work under Assumption \textup{\ref{assuC}}. Assume, moreover, that $\pi\in{\mathcal
R}(\kappa_2)$ with $\kappa_2 >1$. We have
\[
\sup_{a \in(0,1)}\beta(a)<+\infty.
\]
\end{lemm}
\begin{pf}
By Assumption \ref{assuC}, $x \leadsto\pi(x)$ is continuous on $(0,+\infty)$,
hence $\beta(a) = a^{-1}\pi(-{\log a})$ is continuous on $(0,1)$ and it
suffices to show that $\beta(\cdot)$ is bounded in the vicinity of
$0$ and $1$. By assumption, $\pi(x) \lesssim \mathrm{e}^{-\vartheta x}$ for
some $\vartheta\geq1$ near $+\infty$, so $\beta(a) \lesssim
a^{\vartheta-1}$ near the origin and this term remains bounded as
$a\rightarrow0$. By assumption, we also have $\pi\in{\mathcal
R}(\kappa_2)$, so $\pi(x) \lesssim x^{\kappa_2-1}$ near the origin,
therefore $\beta(a) \lesssim(-{\log a})^{\kappa_2-1}$ near $1$ and
this term remains bounded as $a \rightarrow1$ since $\kappa_2 >1$.
\end{pf}

Let $0 < b_\varepsilon\rightarrow0$ as $\varepsilon\rightarrow0$.
For $m>0$, define the class
\[
\widetilde{\mathcal C}_{b_\varepsilon}(m):= \{g\in{\mathcal C}(m),
|{\operatorname{supp} (g ) }|\leq m b_\varepsilon\}.
\]
We have the following extension of Theorem \ref{rate empiricalmeasure}.
\begin{propm} \label{reinforcemenempiricalmeasure}
We work under
Assumptions \textup{\ref{assuA}}, \textup{\ref{assuB}} and \textup{\ref{assuC}}.
Assume that $\pi\in\Pi(\kappa_1)\cap
{\mathcal R}(\kappa_2)$ with $\kappa_1, \kappa_2 > 1$. Then, for
every $1 \leq\mu< \kappa_1+1$,
\[
\sup_{g \in\widetilde{\mathcal C}_{b_\varepsilon}(m)}\mathbb{E}\bigl[
\bigl({\mathcal E}_\varepsilon(g) -{\mathcal E}(g) \bigr)^2 \bigr]=\mathrm{o} \bigl(\varepsilon
^{\mu/(\mu+1)}b_\varepsilon
\bigr) .
\]
\end{propm}
\begin{pf}
We carefully revisit steps 2--4 of the proof of Theorem \ref{rate
empiricalmeasure}, under the additional Assumption \ref{assuC}, and we write
$g(\cdot)=g_\varepsilon(\cdot)$ to emphasize that $g(\cdot)$ may
now depend on the asymptotics.

In step 2,
the right-hand side of (\ref{passagepi}) is now bounded by the
following chain of inequalities:
\begin{eqnarray*}
&&\int_{0}^{-{\log(\varepsilon/ u)}}\mathrm{d}U(s)\int_{-{\log(\varepsilon/
u)}-s}^{+\infty}g_\varepsilon(u \varepsilon^{-1}\mathrm{e}^{-x-s} )^2\pi
(x)\,\mathrm{d}x\\
&&\quad= \int_{0}^{-{\log(\varepsilon/u)}}\mathrm{d}U(s)\int_{0}^{\varepsilon u^{-1}
\mathrm{e}^s}g_\varepsilon(x u \varepsilon^{-1}\mathrm{e}^{-s} )^2\beta(x) \,\mathrm{d}x\\
&&\quad\leq\sup_{a\in(0,1)}\beta(a) u^{-1}\varepsilon\int_{[0,-{\log
(\varepsilon/ u)})}\mathrm{e}^s\,\mathrm{d}U(s)\int_0^1 g_\varepsilon(x )^2\,\mathrm{d}x \lesssim
b_\varepsilon\log(u/\varepsilon),
\end{eqnarray*}
where we have used Lemma \ref{betaborne}, the fact that
$|\operatorname{supp}
(g_\varepsilon)| \lesssim b_\varepsilon$ and $U(s) \lesssim s/c(\pi
)$ again. Therefore,
\[
\mathbb{E}[Y_\varepsilon(v,g)^2 | \lambda_\eta(v)=u ]\lesssim
b_{\varepsilon}\log(\eta/\varepsilon),
\]
hence
\[
\mathbb{E}[M_{\varepsilon,\eta}(g)^2 ] \lesssim b_\varepsilon\eta
\log
(\eta/\varepsilon).
\]
In step 3, we replace $g(\cdot)$ by $g_\varepsilon(\cdot)$ in
${\mathcal E}_\varepsilon(g)$ and ${\mathcal E}(g)$. We first consider
the term $I$. We need to be careful when applying Proposition \ref
{renouvellementSgibnev} because $H_{g_\varepsilon}(z)$ now depends on
$\varepsilon$. By the Cauchy--Schwarz inequality, for $z<0$,
\begin{eqnarray*}
|H_{g_\varepsilon}(-z) |
&\leq& \biggl(\int_{-z}^{+\infty}g_\varepsilon(\mathrm{e}^{-x-z} )^2\pi(x)\,\mathrm{d}x
\biggr)^{1/2} \biggl(\int_{-z}^{+\infty}\pi(x)\,\mathrm{d}x \biggr)^{1/2}\nonumber\\
&\lesssim& \mathrm{e}^{z/2} \biggl(\int_{0}^1g_\varepsilon(y )^2\beta(y\mathrm{e}^z)\,\mathrm{d}y
\biggr)^{1/2}\mathrm{e}^{ \kappa_1 z/2} \lesssim b_\varepsilon^{1/2}\mathrm{e}^{z(1+\kappa
_1)/2},\nonumber
\end{eqnarray*}
again using the fact that $\sup_a \beta(a) \lesssim1$.
We can therefore apply Proposition \ref{renouvellementSgibnev} when $0
< \mu' < (1+\kappa_1)/2$ with rate function $\varphi(z):=\exp(\mu'
z)$, dominating function $r(z):=\mathrm{e}^{-(1+\kappa_1)|z|/2}$, test function
$\psi(z):=b_\varepsilon^{-1/2} 1_{z \leq0}H_{g}(z)$ and $F:=\pi$.
We then obtain, along the same lines as in step~3, for $0 < \mu' <
(1+\kappa_1)/2$, the estimate
\[
\mathbb{E}[I ] \lesssim b_\varepsilon^{1/2} (\varepsilon\eta^{-1}
)^{2\mu
'} \log(1/\eta).
\]
For the term $\mathit{II}$, it suffices to prove that both $\xi_v^{-1}\mathbb
{E}_{\xi
_v} [{\mathcal E}_\varepsilon(g_\varepsilon) ]$ and ${\mathcal
E}(g_\varepsilon)$ are smaller in order than $b_\varepsilon^{1/2}$;
recall (\ref{point important}). For the first term, this follows from
the previous bound on $H_{g_\varepsilon}(z)$ and the representation
(\ref{eqrepresentation}). For ${\mathcal E}(g_\varepsilon)$, since
$\pi\in\Pi(\kappa_1)$ with $\kappa_1>1$, we have, successively,
\begin{eqnarray*}
|{\mathcal E}(g_\varepsilon) | & \leq &\frac{1}{c(\pi)}\int_0^1\frac
{|g_\varepsilon(a)|}{a}\int_{\log(1/a)}^{+\infty}\pi(x)\,\mathrm{d}x \,\mathrm{d}a
\nonumber\\
&\lesssim&\int_0^1 |g_\varepsilon(a)| a^{\kappa_1-1}\,\mathrm{d}a
\lesssim\int_0^1 |g_\varepsilon(a)|\,\mathrm{d}a \lesssim b_\varepsilon.
\end{eqnarray*}
We eventually obtain
\[
\mathbb{E}[B_{\varepsilon,\eta}(g)^2 ]\lesssim b_\varepsilon
\bigl(\varepsilon
^{\mu'}+ (\varepsilon\eta^{-1} )^{2\mu'} \bigr)\log(1/\eta)
\]
for any $0 < \mu' < (1+\kappa_1)/2$.
The trade-off between $M_{\varepsilon,\eta}(g_\varepsilon)$ and
$B_{\varepsilon,\eta}(g_\varepsilon)$ yields the rate
\[
\varepsilon^{\max\{2\mu'/(2\mu'+1),\mu'\}}b_\varepsilon\qquad\mbox
{for any } 0 < \mu' < (1+\kappa_1)/2,
\]
which is of the form $\mathrm{o}(\varepsilon^{\mu/(\mu+1)}b_\varepsilon)$
for any $1 \leq\mu< 1+\kappa_1$, hence the result.
\end{pf}

\textit{Completion of proof of Theorem} \ref{upperparametric}.
By the representation formula (\ref{energieg}), we can write
\[
{\mathcal E}(g_{\gamma_\varepsilon})-m_1(\pi)^{-1}=\frac{1}{m_1(\pi
)}\int_{1-\gamma_{\varepsilon}}^{1}\bigl(f_{\gamma_\varepsilon
}(a)-1\bigr)\beta(a)\,\mathrm{d}a,
\]
where the integral is taken over $[1-\gamma_\varepsilon, 1]$ since
$f_{\gamma_\varepsilon}(a)=1$ on $[0,1-\gamma_\varepsilon]$ and
$\beta(\cdot)$ is a density function with respect to the Lebesgue
measure on $(0,1)$.
We further have
\[
\biggl|\int_{1-\gamma_\varepsilon}^1 \bigl(f_{\gamma_\varepsilon}(a)-1 \bigr)\beta
(a)\,\mathrm{d}a \biggr| \lesssim\int_0^{-{\log(1-\gamma_\varepsilon)}}\pi
(x)\,\mathrm{d}x\lesssim\gamma_{\varepsilon}^{\kappa_2}
\]
since $\|f_\gamma\|_\infty\lesssim1$, by (\ref
{bornitudederiveesf}), $\pi\in{\mathcal R}(\kappa_2)$ and $-{\log
(1-x)} \lesssim x$ for small enough $x\geq0$.
We deduce that
%
\begin{equation} \label{firstestimate}
|{\mathcal E}(g_{\gamma_\varepsilon})-m_1(\pi)^{-1} | \lesssim
\gamma_\varepsilon^{\kappa_2}.
\end{equation}
Next, for some $c>0$, $\gamma_{\varepsilon}g_{\gamma_\varepsilon
}\in\widetilde{\mathcal C}_{\gamma_\varepsilon}(c)$, hence, for any
$0 < \mu< \kappa_1$, Proposition \ref{reinforcemenempiricalmeasure}
entails that
%
\begin{equation} \label{stochpart}\mathbb{E}[ |{\mathcal
E}_\varepsilon
(g_{\gamma_\varepsilon})-{\mathcal E}(g_{\gamma_\varepsilon}) | ]
\lesssim\gamma_\varepsilon^{-1/2}\varepsilon^{\mu/(2\mu+2)}.
\end{equation}
Moreover,
\[
g_{\gamma_\varepsilon}'(a)=-f'_{\gamma_\varepsilon}(a)-af''_{\gamma
_\varepsilon}(a),
\]
hence, by property (\ref{bornitudederiveesf}), we have $\gamma
_\varepsilon^{2}g_{\gamma_\varepsilon} \in{\mathcal C}'(c)$ for
some $c>0$. Applying (\ref{convergencebis}) of Theorem \ref{rate
empiricalmeasure}, we deduce that
%
\begin{equation} \label{additionalerror}
\mathbb{E}[ |{\mathcal E}_\varepsilon(g_{\gamma_\varepsilon
})-{\mathcal
E}_{\varepsilon,\sigma}(g_{\gamma_\varepsilon}) | ] \lesssim\gamma
_\varepsilon^{-2} [ (\sigma\varepsilon^{-1} )^{1/2}+\varepsilon
^{\mu'/4} ]
\end{equation}
for any $0 < \mu' < \kappa_1$.
The specification $\gamma_\varepsilon=\varepsilon^{\mu/(\mu
+1)(2\kappa_2+1)}$ yields the correct rate for (\ref{firstestimate})
and (\ref{stochpart}). The assumption that $\sigma\varepsilon^{-3}$
is bounded ensures that the term $\gamma_\varepsilon^{-2} (\sigma
\varepsilon^{-1} )^{1/2}$ in (\ref{additionalerror}) is
asymptotically negligible since $\kappa_2 \geq1$. Using the fact that
$\kappa_1 \geq4$, the term $\gamma_\varepsilon^{-2}\varepsilon
^{\mu/4}$ also proves negligible by taking $\mu'$ sufficiently close
to $4$. The conclusion readily follows for $\widehat m_{1,\varepsilon}$.

We now turn to higher moment estimators. Thanks to the proof for the
case $k=1$, it suffices to show that
\[
m_1(\pi){\mathcal E}_{\varepsilon,\sigma} (\widetilde g_{\gamma
_\varepsilon} )\rightarrow\int_0^1 \biggl(\log\frac{1}{a} \biggr)^k\beta(a)\,\mathrm{d}a
\]
in probability with the correct rate as $\varepsilon\rightarrow0$.
Note, first, that by representation (\ref{energieg}),
\begin{eqnarray*}
{\mathcal E}(\widetilde g_{\gamma_{\varepsilon}}) &=& \frac{1}{m_1(\pi
)}\int_0^1 h_{\gamma_\varepsilon}(a)\beta(a)\,\mathrm{d}a
\\
&=&\frac{1}{m_1(\pi)} \int_0^1 f_{\gamma_\varepsilon}(1-a) \biggl(\log
\frac{1}{a} \biggr)^k \beta(a)\,\mathrm{d}a,
\end{eqnarray*}
therefore
\[
m_1(\pi){\mathcal E}(\widetilde g_{\gamma_{\varepsilon}})-\int_0^1
\biggl(\log\frac{1}{a} \biggr)^k\beta(a)\,\mathrm{d}a = \int_0^{\gamma_\varepsilon}
\bigl(f_{\gamma_\varepsilon}(1-a)-1 \bigr) \biggl(\log\frac{1}{a} \biggr)^k\beta(a)\,\mathrm{d}a
\]
since $f_{\gamma_\varepsilon}(1-a)=1$ if $a \geq\gamma_\varepsilon
$. It follows that
\begin{eqnarray*}
\biggl|\int_0^{\gamma_\varepsilon} \bigl(f_{\gamma_\varepsilon}(1-a)-1 \bigr)
\biggl(\log\frac{1}{a} \biggr)^k\beta(a)\,\mathrm{d}a \biggr|
&\lesssim&\int_0^{\gamma_\varepsilon} \biggl(\log\frac{1}{a} \biggr)^k\beta(a)\,\mathrm{d}a
\lesssim\int_{\log{1}/{\gamma_\varepsilon}}^{+\infty
}x^k\pi(x)\,\mathrm{d}x \\
&\lesssim& \biggl(\int_{-{\log\gamma_\varepsilon}}^{+\infty} \pi(x)\,\mathrm{d}x
\biggr)^{1-\delta'} \biggl(\int_0^{+\infty} x^{k/\delta'}\pi(x)\,\mathrm{d}x
\biggr)^{\delta'}\\
&\lesssim&\gamma_\varepsilon^{\kappa_1(1-\delta')}
\end{eqnarray*}
for any $0 < \delta' <1$, by H\"older's inequality and where we have
used the fact that $\pi\in\Pi(\kappa_1)$. The second integral in
the last line is finite by Assumption \ref{assuC}. Since the choice of $\delta'$
is free, the choice of $\gamma_\varepsilon$ and the assumption that
$\kappa_1 > \kappa_2$ show that this term is asymptotically
negligible with respect to $(\varepsilon^{\mu/(\mu+1)})^{\kappa
_2/(2\kappa_2+1)}$. Therefore, it suffices to show that
\[
{\mathcal T}_\varepsilon= {\mathcal E}_{\varepsilon, \sigma
}(\widetilde g_{\gamma_{\varepsilon}}) -{\mathcal E}(\widetilde
g_{\gamma_{\varepsilon}})
\]
has order $(\varepsilon^{\mu/(\mu+1)})^{\kappa_2/(2\kappa_2+1)}$.
We split ${\mathcal T}_\varepsilon= {\mathcal T}_{\varepsilon
,1}+{\mathcal T}_{\varepsilon,2}$ with
\[
{\mathcal T}_{\varepsilon,1} = {\mathcal E}_{\varepsilon, \sigma
}(\widetilde g_{\gamma_{\varepsilon}}) -{\mathcal E}_{\varepsilon
}(\widetilde g_{\gamma_{\varepsilon}})
\quad\mbox{and}\quad
{\mathcal T}_{\varepsilon,2} = {\mathcal E}_{\varepsilon}(\widetilde
g_{\gamma_{\varepsilon}}) -{\mathcal E}(\widetilde g_{\gamma
_{\varepsilon}}).
\]
\begin{lemm} \label{controlegprime}
There exists some constant $c>0$, independent of $\varepsilon$, such that:
\begin{itemize}
\item we have $\gamma_\varepsilon^2 \widetilde g_{\gamma_\varepsilon
} \in{\mathcal C}'(c)$;
\item the decomposition
%
\begin{equation} \label{decompositiongtilde}
\widetilde g_{\gamma_\varepsilon}(a) = q_{1,\gamma_\varepsilon
}(a)+q_{2,\gamma_\varepsilon}(a)
\end{equation}
holds, so that for any $0 < \delta' < 1$, we have
$\gamma_\varepsilon^{\delta'}q_{1,\gamma_\varepsilon} \in
\widetilde{\mathcal C}_{\gamma_\varepsilon}(c)$
and
$\gamma_{\varepsilon}^{\delta'} q_{2,\gamma_\varepsilon} \in
{\mathcal C}(c)$.
\end{itemize}
\end{lemm}
\begin{pf}
Tedious but straightforward computations show that
\begin{eqnarray*}
\widetilde g'_{\gamma_\varepsilon} (a) & = & c_{k,1}a^{-1}(\log a
)^{k-2}f_{\gamma_\varepsilon}(1-a)\\
&&{} + [c_{k,2}(\log a)^{k}+c_{k,3}(\log a)^{k-1} ]f'_{\gamma_\varepsilon
}(1-a)+c_{k,4}a(\log a)^{k}f''_{\gamma_\varepsilon}(1-a)
\end{eqnarray*}
with explicit constants $c_{k,1}=(-1)^{k+1}k(k-1)$, $c_{k,2}=(-1)^{k}$,
$c_{k,3}=(-1)^{k}(k+1)k$ and $c_{k,4}=(-1)^{k+1}$. Using property
(\ref{controlefgammaorigine}) of $f_{\gamma_\varepsilon}$, one
readily checks that the four terms multiplied by $\gamma_\varepsilon
^2$ are bounded. For the last term, corresponding to the constant
$c_{k,4}$, the property (\ref{bornitudederiveesf}) of $f_{\gamma
_\varepsilon}$ also shows that this term multiplied by $\gamma
_{\varepsilon}^2$ has the correct order, so $\gamma_\varepsilon^2
\widetilde g_{\varepsilon} \in{\mathcal C}'(c)$ for some $c>0$.

For the second part of the lemma, we have (\ref{decompositiongtilde}) with
\[
q_{1,\gamma_\varepsilon}(a) = (-1)^k a f_{\gamma_\varepsilon}'(1-a)
(\log a )^k
\]
and
\[
q_{2,\gamma_\varepsilon}(a) = (-1)^{k+1}f_{\gamma_\varepsilon
}(1-a)k (\log a )^{k-1}.
\]
By construction of $f_{\gamma_\varepsilon}$, we have $\operatorname
{supp} (q_{1,\gamma_\varepsilon} ) \subset[0,\gamma_\varepsilon]$.
It follows that for any $0 < \delta' <1$ and $a\in(0,1)$, we have
\[
|q_{1,\gamma_\varepsilon}(a) | \leq a^{\delta'} | {\log a} |^k
a^{1-\delta'} |f'_{\gamma_\varepsilon}(1-a) |\lesssim\gamma
_\varepsilon^{1-\delta'} \|f'_{\gamma_\varepsilon}\|_\infty
\lesssim\gamma_\varepsilon^{-\delta'},
\]
where we have used the fact that $\sup_{a\in(0,1)}a^{\delta'} |{\log
a }|^k<+\infty$, the fact that $\operatorname{supp} (q_{1,\gamma
_\varepsilon} ) \subset[0,\gamma_\varepsilon]$ and property (\ref
{bornitudederiveesf}). We conclude that $\gamma_\varepsilon^{\delta
'}q_{1,\gamma_\varepsilon} \in\widetilde{\mathcal C}_{\gamma
_\varepsilon}(c)$ for some $c>0$.

For the term $q_{2,\gamma_\varepsilon}$, we have, for any $a \in
(0,\gamma_\varepsilon]$ and any $0 < \delta' < 1$,
\[
|q_{2,\gamma_\varepsilon}(a) | \leq k a^{\delta'} |{\log a }|^{k-1}
a^{1+\delta-\delta'} \biggl(\biggl(\frac{\gamma}{a}\biggr)^{1+\delta}f_{\gamma
}(1-a) \biggr)
\lesssim1,
\]
where we have again used the fact that $\sup_{a\in(0,1)}a^{\delta'}
|{\log a }|^k<+\infty$ and property (\ref{controlefgammaorigine}).
For $a \geq\gamma_\varepsilon$, we directly have $ |q_{2,\varepsilon
}(a) | \lesssim| {\log\gamma_\varepsilon}|^{k-1}$, which is smaller
in order than $\gamma_\varepsilon^{-\delta'}$ as $\varepsilon
\rightarrow0$.
\end{pf}

The first part of Lemma \ref{controlegprime} enables us to apply
(\ref{convergencebis}) of Theorem \ref{rate empiricalmeasure}: we obtain
\[
\mathbb{E}[|{\mathcal T}_{\varepsilon,1}| ] \lesssim\gamma
_\varepsilon
^{-2} [ (\sigma\varepsilon^{-1} )^{1/2}+\varepsilon^{\mu'/4} ]
\]
for any $0 < \mu' < \kappa_1$ and this term is asymptotically
negligible in the same way as for (\ref{additionalerror}).
The second part of Lemma \ref{controlegprime} enables us to apply
Proposition \ref{reinforcemenempiricalmeasure} to the term
$q_{1,\gamma_\varepsilon}$ and Theorem \ref{rate empiricalmeasure}
to the term $q_{2,\gamma_\varepsilon}$, respectively. It follows that
\begin{eqnarray*}
\mathbb{E}[|{\mathcal T}_{2,\varepsilon}| ]& \leq &\mathbb
{E}[|q_{1,\gamma
_\varepsilon}| ]+\mathbb{E}[|q_{2,\gamma_\varepsilon}| ] \\
& \lesssim &\gamma_\varepsilon^{1/2} \gamma_\varepsilon^{-\delta'}
\varepsilon^{\mu/2(\mu+1)} + \gamma_\varepsilon^{-\delta'}
\varepsilon^{\mu/2(\mu+1)} \lesssim\gamma_\varepsilon^{-\delta
'} \varepsilon^{\mu/2(\mu+1)}.
\end{eqnarray*}
One readily checks that the choice $\delta'<1/2$ shows that this term
is negligible. The proof of Theorem \ref{upperparametric} is thus complete.

\subsection[Proof of Theorem 3]{Proof of Theorem \protect\ref{lowerparametric}}\label{sec54}

Without loss of generality, we consider the homogeneous case with
$\alpha=0$.
We may also assume that $\sigma= 0$ since adding experimental noise to
the observation of the fragments only increases the error bounds.

\textit{Step} 1: \textit{An augmented experiment}. In the binary case, the
dislocation measure $\nu(d\mathbf{s})$ is equivalently expressed via %
a probability measure on $[1/2,1]$ with density function $a \leadsto
\rho(a)$; see (\ref{rho-representation}).

We prove a lower bound in the augmented experiment, where one can
observe all of the sizes $\widetilde X_\varepsilon$ of the fragments
until they become smaller than $\varepsilon$, namely,
\[
\widetilde X_\varepsilon:= \{\xi_u, \xi_{u-} \geq\varepsilon\}
\cup\{\xi_u, u \in{\mathcal U}_\varepsilon\}.
\]
Clearly, taking the infimum over all estimators based on $\widetilde
X_\varepsilon$ instead of $X_{\varepsilon}= X_{\varepsilon, 0}$ only
reduces the lower bound.

For every $u\in{\mathcal U}_\varepsilon$, we have $\xi_{u-} \geq
\varepsilon$. By the conservative Assumption \ref{assuB}, there are at most
$\varepsilon^{-1}$ such $\xi_{u-}$, so $\operatorname{Card}{\mathcal U}_\varepsilon
\leq2\varepsilon^{-1}$. For every node $u \in{\mathcal U}$, the
fragmentation process gives rise to two offspring 
with sizes
$\xi_{u}U$ and $\xi_u(1-U)$, where $U$ is a random variable
independent of $\xi_u$ with density function $\rho(\cdot)$.
Therefore, the process of the sizes of the fragments in the enlarged
experiment can be realized by fewer than 
\[
2\varepsilon^{-1} \biggl(1+\frac{1}{2}+\cdots+\frac{1}{2^{k(\varepsilon
)}} \biggr) \leq\lfloor4\varepsilon^{-1}\rfloor+1=:n(\varepsilon)
\]
independent realizations of the law $\rho(\cdot)$, where
$k(\varepsilon) := \log_2(2/\varepsilon)$, assumed to be an integer
with no loss of generality.

In turn, Theorem \ref{lowerparametric} reduces to proving that
$\varepsilon^{1/2}$ is a lower rate of convergence for estimating
$m_k(\pi)$ based on the observation of an $n(\varepsilon)$-sample of
the law $\rho(\cdot)$. The one-to-one correspondence between $\rho
(\cdot)$ and $\pi(\cdot)$ is given in (\ref{rhopi-correspondence}).

\textit{Step} 2: \textit{Construction of $\pi_\varepsilon$}. We write $\rho_{\pi
}(\cdot)$ to emphasize the dependence on $\pi(\cdot)$. Let
\[
\phi_k(a):=a\log(1/a)^k+(1-a)\log\bigl(1/(1-a) \bigr)^k,\qquad a \in[1/2,1].
\]
%
From (\ref{rhobetta-correspondence}), 
we have
\[
m_k(\pi_0)=\int_{1/2}^1\phi_k(a)\rho_{\pi_0}(a)\,\mathrm{d}a.
\]
Let $0 < \tau< 1$. Choose a function $\psi_k(\cdot
)\dvtx[1/2,1]\rightarrow\mathbb{R}$ such that
\[
\|\psi_k\|_\infty\leq\tau\inf_a \rho_{\pi_0}(a),\qquad \int
_{1/2}^1\psi_k(a)\,\mathrm{d}a=0,\qquad r(k):=\int_{1/2}^1\phi_k(a)\psi_k(a)\,\mathrm{d}a \neq0,
\]
a choice which is obviously possible thanks to Assumption \ref{assuD}. For
$\varepsilon>0$, define
\[
\rho_{\pi_\varepsilon}(a):=\rho_{\pi_0}(a)+\varepsilon^{1/2}\psi
_k(a),\qquad a \in[1/2,1].
\]
(Therefore, (\ref{rhopi-correspondence}) defines $\pi_\varepsilon
(\cdot)$ unambiguously.) By construction, $\rho_{\pi_\varepsilon
}(\cdot)$ is a density function on $[1/2,1]$ and has a corresponding
binary fragmentation with L\'evy measure given by $\pi_\varepsilon
(\cdot)$. Moreover,
\[
m_k(\pi_\varepsilon)=m_k(\pi_0)+r(k)\varepsilon^{1/2}.
\]

\textit{Step} 3: \textit{A two-point lower bound}. The following chain of arguments
is fairly classical. We denote by $\widetilde{\mathbb{P}}_{\pi}$ the
law of
the independent random variables $ (U_i, i=1,\ldots, n(\varepsilon)
)$ with common density $\rho_\pi(\cdot)$ that we use to realize the
augmented experiment.

Let $F_\varepsilon$ be an arbitrary estimator of $m_k(\pi)$ 
based on $\widetilde X_{\varepsilon}$. Put $c:=|r(k)|/2$. We have
\begin{eqnarray*}
&&\max_{\pi\in\{\pi_0,\pi_\varepsilon\}}\widetilde{\mathbb
{P}}_{\pi}
[\varepsilon^{-1/2}|F_\varepsilon-m_k(\pi)|\geq c ]\\
&&\quad\geq \tfrac{1}{2} \bigl(\widetilde{\mathbb{P}}_{\pi_0} [\varepsilon
^{-1/2}|F_\varepsilon-m_k(\pi_0)|\geq c ]+\widetilde{\mathbb{P}}_{\pi
_\varepsilon} [\varepsilon^{-1/2}|F_\varepsilon-m_k(\pi_\varepsilon
)|\geq c ] \bigr)\\
&&\quad\geq \tfrac{1}{2} \widetilde{\mathbb{E}}_{\pi_0} \bigl[1_{\{\varepsilon
^{-1/2}|F_\varepsilon-m_k(\pi_0)|\geq c\}}+1_{\{\varepsilon
^{-1/2}|F_\varepsilon-m_k(\pi_\varepsilon)|\geq c\}} \bigr]-\tfrac
{1}{2}\|\widetilde{\mathbb{P}}_{\pi_0}-\widetilde{\mathbb{P}}_{\pi
_\varepsilon}\|_{\mathrm{TV}},
\end{eqnarray*}
where $\|\cdot\|_{\mathrm{TV}}$ denotes the total variation distance between
probability measures. By the triangle inequality, we have
\[
\varepsilon^{-1/2} \bigl(|F_\varepsilon-m_k(\pi_0)|+|F_\varepsilon
-m_k(\pi_\varepsilon)| \bigr)\geq|r(k)| = 2c,
\]
so one of the two indicators within the expectation above must be equal
to one with full $\widetilde{\mathbb{P}}_{\pi_0}$-probability. Therefore,
\[
\max_{\pi\in\{\pi_0,\pi_\varepsilon\}}\widetilde{\mathbb{P}}_{\pi}
[\varepsilon^{-1/2}|F_\varepsilon-m_k(\pi)|\geq c ] \geq\tfrac
{1}{2}(1-\|\widetilde{\mathbb{P}}_{\pi_0}-\widetilde{\mathbb{P}}_{\pi
_\varepsilon
}\|_{\mathrm{TV}})
\]
and Theorem \ref{lowerparametric} is proved if
%
\begin{equation} \label{TVcontrole}
\limsup_{\varepsilon\rightarrow0}\|\widetilde{\mathbb{P}}_{\pi
_0}-\widetilde{\mathbb{P}}_{\pi_\varepsilon}\|_{\mathrm{TV}}<1.
\end{equation}
By Pinsker's inequality,
$\|\widetilde{\mathbb{P}}_{\pi_0}-\widetilde{\mathbb{P}}_{\pi
_\varepsilon}\|_{\mathrm{TV}}
\leq\frac{\sqrt{2}}{2} (\widetilde{\mathbb{E}}_{\pi_0}
[\log\frac{\mathrm{d}
\widetilde{\mathbb{P}}_{\pi_0}}{\mathrm{d} \widetilde{\mathbb{P}}_{\pi
_\varepsilon}} ] )^{1/2}
$
and
\begin{eqnarray*}
\widetilde{\mathbb{E}}_{\pi_0} \biggl[\log\frac{\mathrm{d} \widetilde{\mathbb
{P}}_{\pi_0}}{\mathrm{d}
\widetilde{\mathbb{P}}_{\pi_\varepsilon}} \biggr]
&=& -\sum_{i=1}^{n(\varepsilon)}\widetilde{\mathbb{E}}_{\pi_0} \biggl[\log
\frac
{\rho_{\pi_\varepsilon}(U_i)}{\rho_{\pi_0}(U_i)} \biggr] \\
&=& -\sum_{i=1}^{n(\varepsilon)}\widetilde{\mathbb{E}}_{\pi_0} \bigl[\log
\bigl(1+\varepsilon^{1/2}\psi_k(U_i)\rho_{\pi_0}(U_i)^{-1} \bigr)-\varepsilon
^{1/2}\psi_k(U_i)\rho_{\pi_0}(U_i)^{-1} \bigr],
\end{eqnarray*}
where we have used the fact that
$\widetilde{\mathbb{E}}_{\pi_0} [\psi_k(U_i)\rho_{\pi_0}(U_i)^{-1}
] = \int
_{1/2}^1\psi_k(a)\,\mathrm{d}a=0$.
We also have that the term $\varepsilon^{1/2}|\psi_k(U_i)\rho_{\pi
_0}(U_i)^{-1}|$ is smaller than $\tau\varepsilon^{1/2}$. Hence, for
small enough $\tau$,
\[
\bigl|-{\log\bigl(1+\varepsilon^{1/2}\psi_k(U_i)\rho_{\pi_0}(U_i)^{-1}
\bigr)+\varepsilon^{1/2}\psi_k(U_i)\rho_{\pi_0}(U_i)^{-1} }\bigr|\leq\tau
^2\varepsilon.
\]
Therefore $\|\widetilde{\mathbb{P}}_{\pi_0}-\widetilde{\mathbb
{P}}_{\pi_\varepsilon
}\|_{\mathrm{TV}} \leq\frac{\sqrt{2}}{2}\tau\varepsilon^{1/2}
n(\varepsilon)^{1/2}$ and this quantity is bounded away from $1$ by
choosing $\tau$ small enough, uniformly in $n$, so (\ref{TVcontrole})
follows. The proof of Theorem \ref{lowerparametric} is thus complete.

\subsection[Proof of Theorem 4]{Proof of Theorem \protect\ref{upperbound}}\label{sec55}

We plan to use the following decomposition:
\[
\widehat\beta(a)-\beta(a) = \widehat m_{1,\varepsilon} {\mathcal
E}_{\varepsilon,\sigma} (h_{a,\varepsilon}
)-\beta(a)=I+\mathit{II}+\mathit{III}+\mathit{IV}
\]
with
\begin{eqnarray*}
I &:=& \widehat m_{1,\varepsilon} \bigl({\mathcal E}_{\varepsilon,\sigma}
(h_{a,\varepsilon}) -{\mathcal E}_\varepsilon(h_{a,\varepsilon} )
\bigr), \\
\mathit{II} & := &\widehat m_{1,\varepsilon} \bigl({\mathcal E}_\varepsilon
(h_{a,\varepsilon} )-{\mathcal E} (h_{a,\varepsilon} ) \bigr), \\
\mathit{III} & := &\bigl(\widehat m_{1,\varepsilon}-m_1(\pi) \bigr) {\mathcal E}
(h_{a,\varepsilon} ),\\
\mathit{IV} & := &m_{1}(\pi) {\mathcal E} (h_{a,\varepsilon} )-\beta(a).
\end{eqnarray*}
Considering \textit{I} and \textit{II}, the term $\widehat m_{1,\varepsilon}$ is bounded
in probability by Theorem \ref{upperparametric}.
By (\ref{convergencebis}) in Theorem~\ref{rate empiricalmeasure},
together with the fact that $\gamma_\varepsilon^{3}\varphi_{\gamma
_\varepsilon,a}' \in{\mathcal C}'(\|\varphi''\|_\infty)$,
we have
%
\begin{equation} \label{estimeI}
\mathbb{E}[ |{\mathcal E}_\varepsilon(h_{a,\varepsilon} )-{\mathcal
E}_{\varepsilon,\sigma} (h_{a,\varepsilon} ) | ] \lesssim\gamma
_\varepsilon^{-3} [ (\sigma\varepsilon^{-1} )^{1/2}+\varepsilon
^{\mu'/4} ]
\end{equation}
for any $0 < \mu' < \kappa_1$.
By construction, we have $\gamma_\varepsilon^{2}\cdot\varphi
'_{\gamma_\varepsilon,a}(\cdot) \in\widetilde{\mathcal C}_{\gamma
_\varepsilon}(\|\varphi'\|_\infty)$.
Therefore, by Proposition \ref{reinforcemenempiricalmeasure},
%
\begin{equation} \label{estimeII}
\mathbb{E}\bigl[ \bigl({\mathcal E}_\varepsilon(h_{a,\varepsilon} )-{\mathcal E}
(h_{a,\varepsilon} ) \bigr)^2 \bigr] \lesssim\gamma_\varepsilon
^{-3}\varepsilon^{\mu/(\mu+1)}.
\end{equation}
Considering \textit{III}, note that for all $a \in(0,1)$, the function $\varphi
_{\gamma_\varepsilon,a}(\cdot)$ has support in $(0,1)$ for
sufficiently small $\varepsilon$ since $\gamma_\varepsilon
\rightarrow0$. Using the representation (\ref{energieg}), we then have
\[
|{\mathcal E} (h_{a,\varepsilon} ) | = \biggl|\frac{1}{m_1(\pi)}\int_0^1
\varphi_{\gamma_\varepsilon,a}(u)\beta(u)\,\mathrm{d}u \biggr| \\
\lesssim m_1(\pi)^{-1} \sup_{u \in(0,1)}\beta(u)
\]
since $\int_0^1 \varphi_{\gamma_\varepsilon, a}(u)\,\mathrm{d}u = \int_0^1
\varphi(u)\,\mathrm{d}u=1$. Recall that $\sup_{u\in(0,1)} \beta(u) \lesssim
1$, by Lemma \ref{betaborne}. By Theorem \ref{upperparametric}, we
conclude that $\mathit{III}^2$ has order
%
\begin{equation} \label{estimeIII}
\varepsilon^{2\mu\kappa_2/(\mu+1)(2\kappa_2+1)}
\end{equation}
in probability. For $\mathit{IV}$, we first note that
$m_{1}(\pi) {\mathcal E} (h_{a,\varepsilon} )=\int_0^1\varphi
_{\gamma_\varepsilon,a}(u)\beta(u)\,\mathrm{d}u$,
hence
\[
\mathit{IV}^{2}= \biggl(\int_0^1\varphi_{\gamma_\varepsilon,a}(u)\beta(u)\,\mathrm{d}u -
\beta(a) \biggr)^2.
\]
The following argument is classical in nonparametric estimation: since
$\beta\in\Sigma(s)$ with $s=n+\{s\}$, where $n$ is a non-negative
integer, by a Taylor expansion up to order $n$ (recall that the number
$N$ of vanishing moments of $\varphi(\cdot)$, recall (\ref
{cancellationproperty}), satisfies $N>s$), we obtain
%
\begin{equation} \label{estimeIV}
\mathit{IV}^2 \lesssim\gamma_\varepsilon^{2s};
\end{equation}
see, for instance, Tsybakov \cite{Tsybakov2}, Proposition 1.2.
Combining (\ref{estimeII}) and (\ref{estimeIV}), we see that the
balance term $\gamma_\varepsilon=\varepsilon^{\mu/(\mu+1)(2s+3)}$
yields the correct rate for $\mathit{II}$ and $\mathit{IV}$. Next, the condition $\kappa
_2 \geq s/3$ ensures that the term (\ref{estimeIII}) also has the
correct order. Finally, the estimate (\ref{estimeI}) proves
asymptotically negligible,
thanks to the assumption that $\sigma\varepsilon^{-3}$ is bounded and
using the fact that $\kappa_1 \geq4$, in the same way as for (\ref
{additionalerror}) in the proof of Theorem \ref{upperparametric}. The
proof of Theorem \ref{upperbound} is thus complete.

\section*{Acknowledgements}
The careful reading, comments and suggestions of two referees helped to
improve considerably an earlier version of this paper. This work was
supported in part by the Agence Nationale de la Recherche, Grant No.
ANR-08-BLAN-0220-01.

\printhistory

\end{document}